\documentclass[a4paper]{amsart}

\usepackage{amsmath}
\usepackage{amsthm} 
\usepackage{mathtools} 
\usepackage{amssymb} 
\usepackage{amsxtra} 
\usepackage{amsfonts} 
\usepackage{txfonts} 
\usepackage{enumerate} 
\usepackage{bbm} 
\usepackage{dsfont} 
\usepackage{bm} 
\usepackage{thm-restate}

\usepackage[T1]{fontenc} 
\usepackage{textcomp} 
\usepackage[utf8]{inputenc}

\newtheoremstyle{break}
{\topsep}{\topsep}
{\itshape}{}
{\bfseries}{}
{\newline}{}
\usepackage{hyperref}

\usepackage{moresize} 
\usepackage{pdfpages} 
\definecolor{mypurple}{HTML}{7f00d4}
\definecolor{mygreen}{HTML}{71a300}
\hypersetup{
	colorlinks=true,
	linkcolor=mypurple,
	filecolor=purple,
	citecolor=mygreen,
	urlcolor=mypurple,
}

\usepackage{marginnote}

\usepackage{graphicx} 
\usepackage{caption} 
\usepackage{subcaption}

\usepackage{pgf,tikz,pgfplots}

\usepackage{tikz-cd}

\usepackage{float} 
\usepackage{turnstile} 
\usepackage{thmtools}
	
	\theoremstyle{break}
	\newtheorem{Teo}{Theorem}[section]
	\newtheorem{mainthm}{Theorem}

	\theoremstyle{definition}
	\newtheorem{Def}[Teo]{Definition}

	\theoremstyle{plain}
	\newtheorem{question}{Question}
	\newtheorem{Lema}[Teo]{Lemma}
	\newtheorem{Prop}[Teo]{Proposition}
	\newtheorem{Coro}[Teo]{Corollary}

	\newenvironment{Dem}{\noindent \vspace{-0.1cm}{\it Proof:}}{\hfill $\square$ \vspace{0.5 cm}} 
	\newenvironment{Not}{\vspace{0.3cm}\noindent \textsc{Notation:}}{\vspace{0.3cm}}{}
	\newenvironment{Remark}{\vspace{0.3cm} \noindent\textsc{Remark.}}{\vspace{0.3cm}}{}
	 
	\newcommand{\concept}[1]{\textbf{#1}}

	\newcommand{\forc}[2]{\dststile{#2}{#1}}
	\newcommand{\tupp}[1]{\left\langle#1\right\rangle}
	\newcommand{\tup}[1]{\langle#1\rangle}
	
	\newcommand{\set}[1]{\{#1\}}

	\newcommand{\RR}{\mathbb{R}}
	\newcommand{\QQ}{\mathbb{Q}}

	\DeclareMathOperator{\OR}{OR}
	
	\DeclareMathOperator{\dom}{dom}
	
	\DeclareMathOperator{\coor}{coor}

	\title[PUCs with no well-ordering of the reals]{Partitions of $\RR^3$ into unit circles with no well-ordering of the reals}
	\author{Azul Fatalini}
	\date{}
	
\begin{document}
		\maketitle
\begin{abstract}
Using a well-ordering on the reals, one can prove there exists a partition of the three-dimensional Euclidean space into unit circles (PUC).  We show that the converse does not hold: there exist models of $\mathsf{ZF}$ without a well-ordering of the reals in which such partition exists. 
Specifically, we prove that the Cohen model has a PUC and construct a model satisfying $\mathsf{DC}$ where this is also the case. 
Furthermore, we present a general framework for constructing similar models for other paradoxical sets, under some conditions of \emph{extendability} and \emph{amalgamation}.
\end{abstract}

{\def\thefootnote{}\footnotetext{The author has been supported by the Deutsche Forschungsgemeinschaft (DFG, German Research Foundation) under the Excellence Strategy EXC 2044–390685587, Mathematics Münster: Dynamics–Geometry–Structure.}}
\section{Introduction}

This paper investigates some paradoxical sets of reals and study their interaction with the Axiom of Choice.
Informally, paradoxical sets are subsets of $\mathbb{R}^n$ that can be constructed using the Axiom of Choice. 
Their existence can be counter-intuitive at first sight: for example, the well-known examples of a non-measurable set by Vitali \cite{Vitali1905} and the partition given by the Banach--Tarski paradox \cite{Banach1924}. 
Although there are many examples of paradoxical sets (see, for example, \cite{Komjath1993}), much remains unknown about these objects. 

In this work we will focus primarily on a specific paradoxical set: 
A partition of $\mathbb{R}^3$ into unit circles (\concept{PUC}).
The known proofs of existence of this object rely on a transfinite induction on a well-order of the reals \cite{Conway1964}.

For any notion of paradoxical set there are natural questions to ask about its properties. 
For example, whether we can have a paradoxical set that is Borel, measurable, meager, etc.
In Subsection \ref{subsection: PUC lit review}, we provide a literature review of PUCs. Many questions about this object remain unanswered, such as the following. 

\begin{question}\label{question: borel}
	Can a $\mathrm{PUC}$ be Borel? \emph{\cite{HamkinsMO}} 
\end{question}

There is another paradoxical set called two-point set, which has a similar flavor and has been studied much more. The proof of its existence also relies on a well-ordering of the reals \cite{Mazurkiewicz1914}.
After extensive efforts to determine whether a two-point set can be Borel \cite{Baston1972, Mauldin1998, Miller1989}, another approach to this question was needed. 
Recent work has shifted focus to studying these objects from a set-theoretical perspective \cite{Beriashvili2022, Brendle2018, Larson2020, Schilhan2022}. 
In that direction, the main question considered through this work is the following.

\begin{question}\label{question: recover choice}
	\textbf{Can we recover some weakening of the Axiom of Choice from the existence of a particular paradoxical set?}
\end{question}

There have been a number of recent results on this topic. 
Larson and Zapletal \cite{Larson2020} developed a broad technique that deals with a similar type of problems, under some large cardinals assumptions.
There has also been some progress using symmetric extensions of models of set theory. 
In particular, Schilhan \cite{Schilhan2022} gives a partial negative answer to Question \ref{question: recover choice} for sets of reals that can be defined as \emph{maximal independent sets}.

This paper takes a different approach from the aforementioned lines of work.
First, we do not have any large cardinal assumption. 
Second, the objects we considered (PUCs and Mazurkiewicz sets) cannot be defined as maximal independent sets. This is usually the case for partitions. 
Instead, the direction of this work follows the lines of other authors (for example, \cite{Beriashvili2022, Beriashvili2018, Brendle2018, Kanovei2020, Schindler2018}).

The contribution of this article is giving negative answers for different versions of Question \ref{question: recover choice}, changing the particular paradoxical set and the weakening of the Axiom of Choice considered. 
The choice principles we will examine are: the existence of a well-ordering of the reals ($\mathsf{WO}(\RR)$), the Principle of Dependent Choices ($\mathsf{DC}$), and Countable Choice ($\mathsf{AC}_\omega$).

For models of $\mathsf{ZF}+\mathsf{DC}$, we develop a framework in Section \ref{section: gral setup} to construct models of that form with a specific paradoxical set given by a definition $\psi$. \\

\begin{mainthm}[see Theorem \ref{th: partition no wo}]
	Let $V$ be a model of $\mathsf{ZFC}$. 
	Let $\mathbf{Q}$ be the finite support product of $\omega_1$-many copies of Cohen forcing, and let $g$ be a $\mathbf{Q}$-generic filter over $V$.
	Let $\mathbf{P}$ be a forcing notion over $V[g]$ that adds a real partition, let $h$ be a $\mathbf{P}$-generic filter over $V[g]$, and let $\mathcal{P}=\cup h$.
	
	If $\mathbf{P}$ is $\sigma$-closed and satisfies extendability and amalgamation, then 
	\[ L(\RR, \mathcal{P})^{V[g,h]} \models \mathsf{ZF}+\mathsf{DC}+\neg \mathsf{WO}(\RR) + \psi(\mathcal{P}). \] 
\end{mainthm}
 
This method can be used to recover known results about Hamel bases and Mazurkiewicz sets (Subsections \ref{subsection: hamel bases} and \ref{subsection: maz sets}), and it can be used for new applications. 
In particular, we apply it to PUCs in Section \ref{section: PUC}, and  obtain the following result.\\

\begin{mainthm}[see Theorem \ref{th: PUC - no wo}]
	There is a model of
	\[  \mathsf{ZF}+\mathsf{DC}+\neg \mathsf{WO}(\RR)+ \text{there is a partition of $\RR^3$ into unit circles}. \]
\end{mainthm}

Since we are considering choiceless models, it is of interest to know what happens in the first example of such: the \emph{Cohen model} (also called \emph{Cohen-Halpern-Lévy model}), a model which does not satisfy Countable Choice.
It is known that in this model there are many examples of paradoxical sets while the Axiom of Choice fails dramatically (see Theorem \ref{th: prereq paradoxical sets in H}).
By refining the definition of \emph{amalgamation} to $\emph{($<\omega$)-amalgamation}$ and with the help of recent results on the transcendence degree of certain set-theoretical subfields \cite{Fatalini2024}, we can prove that there is a PUC in the Cohen model as well. That is the goal of Section \ref{section: PUC in the Cohen model}.\\

\begin{mainthm} [see Theorem \ref{th: PUC in Cohen model}]
There is a partition of $\RR^3$ into unit circles in the Cohen model.
\end{mainthm}

\subsection{Overview of PUCs}\label{subsection: PUC lit review}

In 1964, Conway and Croft analyzed the problem of covering $S^n$ or $\RR^n$ with open/closed/half-closed arcs and segments respectively, of the same length \cite{Conway1964}.
They answered many of these questions and provided several explicit such partitions (in $\mathsf{ZF}$). 
However they could not find an explicit solution to the problem of partitioning $S^n$ into closed arcs of the same length.
They developed a more general theorem \cite[Appendix]{Conway1964} that could be applied to all of these problems for dimension $n\geq 3$, but it used the Axiom of Choice. 
A corollary of this theorem is the existence of a partition of $\RR^3$ into unit circles (here as Theorem \ref{th: zfc PUC}).
Jonsson \cite{Jonsson1998} attributes the latter result also to Kharazishvili, who seems to have proven it in 1985, but the present author could not recover this reference. 

There is no trace that somebody looked into that or similar objects (after all, the result only appears in the last sentence of the appendix in the paper of Conway and Croft), until Szulkin in 1983 showed in a surprisingly simple way that it is possible to partition $\RR^3$ in circles without the Axiom of Choice, but dropping the requirement of the circles to have the same radius \cite{Szulkin1983}.
Furthermore, it is easy to verify that $\RR^2$ cannot be partitioned into circles (not even in Jordan curves).
To the best of the author's knowledge, there are no more results regarding unit circles that did not use the Axiom of Choice until very recently: now we know that there is an open set of $\RR^3$ for which there is an explicit partition of unit circles \cite[Example 3.1]{Asimov2024}.
However, it is still open whether it is possible in $\mathsf{ZF}$ to prove the existence of a partition of the full three-dimensional euclidean space into unit circles.

In 1985, Bankston and Fox 
tried to expand Theorem \ref{th: zfc PUC} (but topologically) to higher dimensions of the euclidean space as well as of the spheres that are used to tile it. 
For similar reasons to the case $n=1$, $S^n$ cannot partition $\RR^{n+1}$ for any $n$, not even allowing the tiles to be \emph{topological} copies of $S^n$ \cite[Theorem 2.3]{Bankston1979}.
But Bankston and Fox proved (in $\mathsf{ZF}$) that $\RR^{n+2}$ (and therefore any bigger dimension) can be partitioned into topological copies of $S^n$ for all $n<\omega$ \cite{Bankston1985}.
Additionally, the proof of Theorem \ref{th: zfc PUC} can be generalized to prove (in $\mathsf{ZFC}$) that $\RR^{2n+1}$ can be partitioned into \emph{isomorphic} copies of $S^n$ \cite[Theorem 2.5]{Bankston1979}.
To the best of the author's knowledge, the question of whether $\RR^m$ can be partitioned into isomorphic copies of $S^n$ is open for $n+2\leq m < 2n+1$ and $n\geq 2$, with or without the Axiom of Choice, and the same is open even when allowing different radii. 
As the simplest example, it is not known whether $\RR^4$ can be partitioned into two dimensional spheres \cite[Question 3.1.iv]{Bankston1979}.

A natural question arises: into what types of pieces can $\RR^3$ be partitioned? 
For example, $\RR^3$ can be partitioned in: letters T \cite{Rosen2022}, rhombi with edge length 1 \cite[Theorem 3]{Wilker1989} (not known for filled squares), \emph{unlinked} circles of the same radius \cite[Theorem 2.3]{Jonsson1998}, unlinked circles where every positive real number appears exactly once as a radius \cite[Theorem 2.1]{Jonsson1998}, and even any family of cardinality $\mathfrak{c}$ of real analytic curves \cite[Theorem 3.1]{Jonsson1998}.
However, it is not true that $\RR^3$ can be partitioned into isometric copies of any fixed Jordan curve.
A nice overview of (subsets of) $\RR^n$ that can be partitioned into (subsets of) $\RR^m$ is displayed by Jonsson and Wästlund \cite[Section 4]{Jonsson1998}.

In yet another direction, Kharazishvili explores the concept of $k$-homogeneous covering, i.e., each point of the euclidean space considered is covered by exactly $k$ tiles. 
For $k=1$, these are just partitions. 
While $\RR^2$ cannot be partitioned by copies of $S^1$, there exists a simple $2$-homogeneous covering of $\RR^2$ in circles of the same radius \cite[Example 2]{Kharazishvili2010}.
We refer the reader to \cite{Kharazishvili2004, Kharazishvili2010} for more examples.

On the side of negative results, Cobb proved in 1995 that
$\mathbb{R}^3$ cannot be continuously decomposed into circles \cite{Cobb1996}.
Continuously, in this case, means that any sequence of points converging to some point $x$ induces a convergence of the radii, planes, and centers of the circles associated with them in the partition, and they converge to those of the circle that passes through $x$. 

Different from other examples of paradoxical sets, we do not know much  about PUCs. 
The contribution of this article is to give similar results that were known for Hamel basis and Mazurkiewicz sets but for the case of PUCs. 
Hamkins asked whether a partition of $\RR^3$ into unit circles can be Borel \cite{HamkinsMO}.
Similarly to the Mazurkiewicz case (for which this question is also open), we will show that in case we find a partition into unit circles that is analytic, then it is Borel (see Proposition \ref{prop: PUC if analytic then borel}). 
Using the strategy of Miller \cite{Miller1989} for obtaining coanalytic Hamel bases and Mazurkiewicz sets under the assumption of $V=L$, which was later generalized by Vidnyánszky \cite[Theorem 3.4]{Vidnyanszky2012}, it can also be shown that if $V=L$, then there is a coanalytic PUC (see Proposition \ref{prop: coanalytic PUC}).
In terms of our guiding question (Question \ref{question: recover choice}), we will exhibit a model of $\mathsf{ZF}+\mathsf{DC}+\neg\mathsf{WO}(\RR)$ with a partition of $\RR^3$ into unit circles (Theorem \ref{th: PUC - no wo}) by applying the methods of Section \ref{section: gral setup}.
Moreover, we will show in Theorem \ref{th: PUC in Cohen model} that the Cohen model has a PUC, so we cannot recover countable choice from the existence of this paradoxical set. 

\subsection{More on paradoxical sets in choiceless models}\label{subsection: lit review paradoxical sets}

In this subsection we give an overview of some facts about other paradoxical sets, Hamel bases and Mazurkiewicz sets, to later compare them and its properties with the case of PUCs (see section  \ref{subsection: properties PUCs}).

\subsubsection{Hamel bases} \label{section: hamel basis lit review}

In $\mathsf{ZF}$, the existence of a Hamel basis (a basis of $\RR$ as a $\QQ$-vector space) implies the existence of a Vitali set, which is the standard example of a subset of $\RR$ which is not Lebesgue-measurable.
However, a Hamel basis itself can be measurable, and every measurable Hamel basis has null measure \cite[Theorem I]{Sierpinski1920}. 
Yet there is no Hamel basis that is Borel \cite[Theorem 2]{Sierpinski1920}.
Actually, it is known that there is no analytic Hamel basis \cite[Theorem 9]{Jones1942}.
The next best possibility is being coanalytic, which is consistent and follows from $V=L$ \cite[Theorem 9.26]{Miller1989}.

The existence of a Hamel basis can be proven from $\mathsf{ZF}+\mathsf{WO}(\RR)$, and it is of course a corollary of the existence of bases for every vector space, which is equivalent to $\mathsf{AC}$.
With Question \ref{question: recover choice} in mind,  a natural question is whether one can recover $\mathsf{WO}(\RR)$ from the existence of a Hamel basis. This was asked by Pincus and  Prikry \cite{Pincus1975}.
It turns out that one cannot recover even Countable Choice, since in the Cohen-Halpern-Lévy model (Definition \ref{def: prereq cohen model}) there is a Hamel basis \cite[Theorem 2.1]{Beriashvili2018}.
If one is interested in having $\mathsf{DC}$ in the model, it is also possible: it is consistent that there is a model of $\mathsf{ZF}+\mathsf{DC}+\neg \mathsf{WO}(\RR)$ with a Hamel basis \cite[Theorem 1.1]{Schindler2018}. 
We will recover this result (see \ref{coro: hamel basis - no wo}) using the results on Section \ref{section: gral setup}.
Moreover, in a follow-up article \cite[Theorem 5.4]{Brendle2018} the authors show that there is a model $M$ of $\mathsf{ZF}+\mathsf{DC}+\neg \mathsf{WO}(\RR)$ in which there is a Hamel basis (moreover, a Burstin basis, i.e. a Hamel basis which has nonempty intersection with every perfect set) and several other paradoxical sets (Luzin, Sierpi\'{n}ski, and of course Vitali).
Using the same methods from Section \ref{section: gral setup} we can also recover the existence of a Hamel basis in the model $M$ (see Corollary \ref{coro: hamel basis - sacks forcing}). 
Furthermore, using an inaccessible cardinal Larson and Zapletal produced a model of $\mathsf{ZF}+\mathsf{DC}+$ in which there is a Hamel basis but no non-principal ultrafilter on $\omega$ \cite[Corollary 12.2.10]{Larson2020}.

\subsubsection{Mazurkiewicz sets}

Using the Axiom of Choice in the form of a transfinite induction on the cardinality of the continuum, one can prove that there is a Mazurkiewicz set (also called two-point set), namely, a subset of the plane that intersects every line in exactly 2 points. 
So only $\mathsf{ZF}+\mathsf{WO}(\RR)$ is needed to carry out the construction.
By a similar proof, one can construct $n$-point subsets of $\RR^2$ for $n<\omega$.
Even more, if for any line $l$ of the plane we assign a cardinal $\alpha_l$ such that $2\leq \alpha_l \leq 2^{\aleph_0}$, then there exists a set of points in the plane that intersects every line $l$ in precisely $\alpha_l$ points \cite{Sierpinski1953}.
Mauldin raised the question of which are the conditions for which we can have a Borel set that meets every line $l$ in exactly $2\leq \alpha_l<\omega$ points. The function $l\mapsto \alpha_l$ needs to be Borel \cite[Theorem 12]{Mauldin1998}, but what else can we say?
There is a simple example of an $\aleph_0$-point set which is $F_\sigma$: the union of all the circles centered in the origin and with integer radii. 
Nevertheless, the case of $\alpha_l$ being a fix natural number $n$ for any $n\geq2$ is still open. 

Unlike the case of Hamel bases where any linearly independent subset of $\RR$ can be extended to a Hamel basis, deciding whether a partial Mazurkiewicz set (a subset of $\RR^2$ without three points on a line) can be extended to a (full) Mazurkiewicz set is very hard. 
The usual proof for the existence of Mazurkiewicz sets shows that any partial two-point set of cardinality strictly less than continuum can be extended to a full two-point set. 
However, there are \emph{small} partial Mazurkiewicz sets of cardinality $\mathfrak{c}$ that cannot be extended to a full Mazurkiewicz set. 
The simplest example of this is a circle. 
A circle has cardinality $\mathfrak{c}$, is a partial two-point set, but it is easy to see that it cannot be extended to a (full) two-point set.
More on this topic is studied by Dijkstra, Kunen and van Mill \cite{Dijkstra1998, Dijkstra1999}.

Going back to the question of whether a Mazurkiewicz set can be Borel, Baston and Bostock discard the simplest case: a two-point set cannot be $F_\sigma$ \cite[Theorem 3]{Baston1972}. 
Similar to the case of Hamel bases, there are measurable and non-measurable Mazurkiewicz sets, and every measurable Mazurkiewicz set has measure zero \cite[II.10.21]{Gelbaum2003}. 
A Mazurkiewicz set must have topological dimension zero \cite[Theorem 2]{Kulesza1992}, which answers a question of Mauldin \cite[1069 Problem 2.3]{Mauldin1990}.
The question whether $3$-point sets have this property seems open, while $n$-point sets for $n\geq4$ may be one-dimensional (they could contain a circle!).

It is also known that if an $n$-point set is analytic then it is Borel \cite[Section 7]{Miller1989}. 
As in the case of Hamel bases, if $V=L$ then there is a coanalytic Mazurkiewicz set \cite[Theorem 7.21]{Miller1989}, and the same proof shows that the same holds for $n$-point sets. \\

Attending to the main question of this paper (Question \ref{question: recover choice}), it is also true that one cannot recover a $\mathsf{WO}(\RR)$ from the existence of a Mazurkiewicz set, as it was the case for Hamel bases. 
The first to show this was Miller \cite[Theorem 5]{Miller2008}.
As in the case of Hamel bases, one cannot even recover countable choice, since the Cohen-Halpern-Lévy model $H$ contains a Mazurkiewicz set \cite[Corollary 0.3]{Beriashvili2022}.
The strategy for proving that $H$ has a Hamel basis \cite{Beriashvili2018} and a partition of unit circles (Theorem \ref{th: PUC in Cohen model}) is proving that the object satisfies \emph{$(<\omega)$-amalgamation} (Definition \ref{def: omega amalagamation}).
In the case of Mazurkiewicz sets, this strategy does not seem to work. 
Nevertheless, Beriashvili and Schindler gave a criteria for a model to have a Mazurkiewicz set by exploiting a geometrical construction of Chad, Knight and Suabedissen \cite[Lemma 4.1]{Chad2009}, which was also used in the construction of the model by Miller.
Furthermore, there is a model of $\mathsf{ZF}+\mathsf{DC}+\neg\mathsf{WO}(\RR)$ with a Mazurkiewicz set \cite{Beriashvili2019}.
We can recover this result by using the methods in Section \ref{section: gral setup}, and this is shown in Subsection \ref{subsection: maz sets}.

\section{Preliminaries}
	 We will assume the reader is familiar with the basics of forcing, but we will explicitly state some definitions and properties that will be used often along this text, and fix some notation.
	 Here, a forcing notion $\mathbf{P}=(\mathbf{P}, \leq_{\mathbf{P}}, \mathbbm{1}_{\mathbf{P}})$ is a partially ordered set with a largest element $\mathbbm{1}_\mathbf{P}$.
	 
	 \begin{Def}\label{def: prereq cohen forcing and C(X)}
	 	Let $\mathbf{C}={^{<\omega}}\omega $ be the forcing given by 
	 	\[	\mathbf{C} = \{p\colon \omega \to \omega \mid  p \text{ is a partial function with } \dom(p)<\omega \},\]
	 	ordered by reverse inclusion and $\mathbbm{1}_{\mathbf{C}}=\emptyset$.
	 	We call this forcing \concept{Cohen forcing}\footnote{Notice that here we denote this forcing by $\mathbf{C}$ and not $\mathbb{C}$ as frequently seen in the literature. This is to avoid confusion with the complex plane $\mathbb{C}$.}. 
	 	
	 	For a set of ordinals $X$, we write $\mathbf{C}(X)$ for the finite support product of $X$-many copies of $\mathbf{C}$. 
	 	Namely, 
	 	
	 	\[\mathbf{C}(X)= \left\{p \in \Pi_{\alpha\in X} \mathbf{C} : \left| \{ \alpha\in X : p(\alpha)\neq \emptyset \} \right| <\omega  \right\}, \]
	 	ordered coordinatewise.
	 \end{Def}
	 
	 If $g\subseteq \mathbf{C}$ is a generic filter over a model $V$, $\cup g \in {^\omega}\omega \cap V[g]$ is a real usually called the \concept{Cohen real} added by $g$.
	 If $g\subseteq \mathbf{C}(X)$ is a generic filter over $V$, then 
	 $\cup(g\restriction \{\alpha\})$ is also a real for each $\alpha \in X$.
	 We will often mix up the generic filters with the reals added by them for these forcings.
	 
	 We will use several times a nice fact of $\mathbf{C}(\omega_1)$ that establishes that any real in a forcing extension by this forcing is in the model produced by some (strict) initial segment of the generic. 
	 
	 \begin{Lema}\label{lemma: prereq any real is in an inital segment of C(omega1)}
	 	Let $g$ be $\mathbf{C}(\omega_1)$-generic over $V$ and let $r\in {^\omega}\omega \cap V[g]$. 
	 	Then there is $\alpha<\omega_1$ such that $r\in V[g\restriction \alpha]$.
	 \end{Lema}
	 
	 \begin{Dem}
	 	Identify $r$ as an element of ${^\omega}2$.
	 	Let $\tau$ be a name for $r$. 
	 	We can assume $\tau$ is of the form $\{ (\check{n}, p) \mid p \in A_n\}$, 
	 	where $A_n$ is an antichain for every $n<\omega$. 
	 	It is a well-known fact that $\mathbf{C}(\omega_1)$ is ccc, namely, every antichain on $\mathbf{C}(\omega_1)$ is countable.
	 	Therefore each $A_n$ is countable. 
	 	For each $n$, let $\alpha_n$ be the supremum of the supports of conditions in $A_n$. 
	 	Since $A_n$ is countable, $\alpha_n<\omega_1$. 
	 	Let $\alpha$ be the supremum of $\{\alpha_n\}_{n<\omega}$.
	 	Again, $\alpha<\omega_1$.
	 	Then $\tau$ is also a name in $V^{\mathbf{C}(\alpha)}$ and therefore $r=\tau_g=\tau_{g\restriction\alpha}\in V[g\restriction \alpha]$.
	 \end{Dem}
	 
	 Notice that $g\restriction\alpha$ is $\mathbf{C}(\alpha)$-generic over $V$; so the notation $V[g\restriction \alpha]$ makes sense.
	 Also, observe that $\mathbf{C}(\omega_1)$ is essentially the same (there is a natural isomorphism of partial orders) as $\mathbf{C}(\alpha)\times \mathbf{C}(\omega_1 \backslash \alpha)$ for each $\alpha<\omega_1$.
	 Moreover, if $X$  has cardinality $\aleph_1$ in $V$, then $\mathbf{C}(X)$ is isomorphic to $\mathbf{C}(\omega_1)$.\\
	 
	 There is another way in which two forcing relations can be related, which is called \emph{forcing equivalence} and is denoted by $\cong$.
	 The following theorem gives us a nice characterization of Cohen forcing up to forcing equivalence.\\
	 
	 \begin{Teo}[{\cite[Theorem 1 Section 4.5]{Grigorieff1975}}] \label{th: prereq cohen is the only countable forcing}
	 	Let $\mathbf{P}$ be a separative, countable and atomless poset.  
	 	Then $\mathbf{P}$ contains a dense subset isomorphic to $\mathbf{C}$.
	 	In other words, Cohen forcing $\mathbf{C}$ is the only countable atomless forcing (modulo forcing equivalence).                                                                                                                                                                             
	 \end{Teo}
	 
	 \begin{Remark}
	 	If $X$ is a set of ordinals that is countable in $V$, then $\mathbf{C}\cong \mathbf{C}(X)$.
	 \end{Remark}
	 
	 It is a well known fact that any forcing is forcing equivalent to a separative forcing \cite[Lemma 14.11]{Jech2003}, so we do not need to check whether $\mathbf{P}$ is separative in Theorem \ref{th: prereq cohen is the only countable forcing}.
	 
	 \begin{Not}
	 	If $M$ and $N$ are two models of $\mathsf{ZF}$, we write $N \hookrightarrow M$ to denote that $M$ is a forcing extension of $N$. 
	 	We write $N \xhookrightarrow{\mathbf{P}} M$ if $N$ is a \concept{$\mathbf{P}$-ground} of $M$, i.e., $M=N[g]$ for some $\mathbf{P}$-generic filter $g$ over $N$.
	 \end{Not}
	 
	 The following is a result that we will use several times. \\
	 
	 \begin{Teo}[The Solovay basis result {\cite[Theorem~2 Section~2.14]{Grigorieff1975}}] 
	 	\label{th: prereq solovay}
	 	Let $M$ be a model of $\mathsf{ZF}$, $\mathbf{P}\in M$ be a forcing notion and let $g$ be a $\mathbf{P}$- generic filter over $M$. 
	 	If $a\in M[g]$ and $a\subseteq M$, then 
	 	\[M \hookrightarrow M[a] \hookrightarrow M[g]. \]
	 	Moreover, the first forcing is given by a complete subalgebra of the completion boolean algebra of $\mathbf{P}$ and the second is a forcing given by the quotient $\mathcal{B}/H$ where $H$ is a generic filter of the first forcing.
	 \end{Teo}
	 
	 Now we will apply Theorem \ref{th: prereq solovay} to our favorite forcings $\mathbf{C}$ and $\mathbf{C}(\omega_1)$. 
	 
	 \begin{Lema}\label{lemma: prereq any real in cohen ext is cohen}
	 	Let $g$ be a $\mathbf{C}$-generic filter over $V$ and let $r\in {^\omega}\omega \cap V[g]$. 
	 	Then 
	 	\[V \xhookrightarrow{a} V[r] \xhookrightarrow{b} V[g],\]
	 	where $a\cong \mathbf{C}$ if $V[r]\neq V$, and $b\cong \mathbf{C}$ if $V[r]\neq V[g]$. 
	 \end{Lema}
	 
	 \begin{Lema}\label{lemma: prereq any real in Q extension is Q ground}
	 	Let $g$ be a $\mathbf{C}(\omega_1)$-generic filter over $V$ and let $r\in {^\omega}\omega \cap V[g]$. 
	 	Then 
	 	\[V \xhookrightarrow{a} V[r] \xhookrightarrow{b} V[g],\]
	 	where $a\cong \mathbf{C}$ if $V[r]\neq V$ and $b\cong \mathbf{C}(\omega_1)$. 
	 \end{Lema}
	 
	 \begin{Dem}
	 	Let $\alpha$ be such that $r\in V[g\restriction \alpha]$ (see Lemma \ref{lemma: prereq any real is in an inital segment of C(omega1)}).
	 	Let $x$ be the Cohen real such that $V[x]=V[g\restriction \alpha]$ (see the remark below Theorem \ref{th: prereq cohen is the only countable forcing}).
	 	By Theorem \ref{th: prereq solovay},  
	 	\[V\xhookrightarrow{a} V[r]  \xhookrightarrow{c} V[x]=V[g\restriction \alpha] \xhookrightarrow{\mathbf{C}(\omega_1\backslash \alpha)} V[g].\]
	 	Moreover, applying Theorem \ref{lemma: prereq any real in cohen ext is cohen}, 
	 	$a$ and $c$ are each forcing equivalent to Cohen forcing or a trivial forcing.
	 	Then $b = \mathbf{C}(\omega_1 \backslash \alpha)$ or $b\cong \mathbf{C}\times \mathbf{C}(\omega_1 \backslash \alpha)$. 
	 	In any case, $b\cong \mathbf{C}(\omega_1)$.
	 \end{Dem}

There is yet another well-known concept which is \emph{mutual genericity}.
	 For example, if $g$ is $\mathbf{C}(\omega_1)$-generic over $M$, then $g\restriction\alpha\subseteq\mathbf{C}(\alpha)$ and $g\restriction(\omega_1\backslash\alpha)\subseteq\mathbf{C}(\omega_1 \backslash \alpha)$ are mutually generic for any $\alpha<\omega_1$. 
	 One property of mutual genericity that we will use very often is the following: if $g$ and $h$ are mutually generic over a model $M$, then $M[g]\cap M[h]=M$.
	 
	 \begin{Def}\label{def: prereq split of a cohen real}
	 	Let $x$ be a function $x\colon \omega\to \omega$.
	 	We can \emph{split} $x$ in two reals $x_0$, $x_1$ such that $x=x_0\oplus x_1$, where $\oplus$ is the operation of alternating digits from each of the reals, namely, $x(2n)=x_0(n)$ and $x(2n+1)=x_1(n)$ for all $n<\omega$.
	 	If $s$ is a finite initial segment of $x$, we say we \concept{split $x$ according to $s$} in two reals $x_0$ and $x_1$ iff $s$ is an initial segment of both  $x_0$ and $x_1$ and $x\backslash s = (x_0 \backslash s)\oplus (x_1\backslash s)$. 
	 \end{Def}
	 
	 \begin{Remark}
	 	Let $x$ be a Cohen real over a model $M$. 
	 	In $M[x]$, let $x_0, x_1$ be the split of $x$ according to $s$. 
	 	Then $x_0$ and $x_1$ are mutually generic Cohen reals over $M$.
	 \end{Remark}
	 
	 Our favorite forcings $\mathbf{C}$ and $\mathbf{C}(\omega_1)$ satisfy a nice property that will be useful for our purposes, namely, they are \emph{homogeneous}.
	 
	 \begin{Def} \label{def: prereq homogeneous} 
	 	Let $\mathbf{P}$ be a poset. 
	 	We call $\mathbf{P}$ \concept{ homogeneous} iff for all $p,q\in \mathbf{P}$, there is a dense homomorphism $\pi$ from $\mathbf{P}$ to itself such that $\pi(p)=q$. 
	 \end{Def}
	 
	 \begin{Lema}[{\cite[Lemma 6.53]{Schindler2014}}] \label{lemma: prereq C and C(alpha) homogeneous}
	 	$\mathbf{C}$ is homogeneous. 
	 	If $\alpha$ is an ordinal, then $\mathbf{C}(\alpha)$ is homogeneous.
	 \end{Lema}
	 
	 The only consequence of being homogeneous that will be needed is stated in the following lemma. 
	 
	 \begin{Lema}[{\cite[Lemma 6.61]{Schindler2014}}] \label{lemma: prereq P homogeneous implies 1 forces}
	 	Let $M$ be a transitive model of $\mathsf{ZFC}$, let $\mathbf{P}\in M$ be a homogeneous forcing notion. 
	 	Let $\phi$ be a formula and $x_0, \dots, x_{n-1} \in M$. 
	 	Then
	 	\[\text{$\mathbbm{1}$} \forc{\mathbf{P}}{M} \phi (\check{x}_0, \dots, \check{x}_{n-1}), \text{ or } \text{$\mathbbm{1}$} \forc{\mathbf{P}}{M} \neg\phi (\check{x}_0, \dots, \check{x}_{n-1}).\]
	 \end{Lema}

	Finally, in Section \ref{section: PUC in the Cohen model} we will prove that there is a PUC in the classical first example of a model in which the Axiom of Choice fails: the First Cohen model, also called the Cohen–Halpern–Lévy model.
	 
	 \begin{Def}\label{def: prereq cohen model}
	 	Let $g$ be a $\mathbf{C}(\omega)$-generic filter over $L$. 
	 	Let us write $A = \{c_n : n < \omega\}$ for the set of Cohen reals added by $g$, i.e., $c_n=\cup (g\restriction\{n\})$ for $n<\omega$. 
	 	The model 
	 	\[H= \mathsf{HOD}_A^{L[g]}\]	
	 	of all sets which are hereditarily ordinal definable inside $L[g]$ from parameters in $A \cup \{A\}$ is called \concept{the Cohen–Halpern–Lévy model}. 
	 \end{Def}
	 
	 We will use this model in Section \ref{section: PUC in the Cohen model} so we will describe some of its properties. 
	 It was introduced by Cohen \cite[pp. 136--141]{Cohen1966}, and explored later in a different presentation by Halpern and Lévy \cite{Halpern1974}.\\
	 
	 \begin{Teo}[{\cite[pp. 136--141]{Cohen1966}}]\label{th: in H A has no countable subset}
	 	In the Cohen–Halpern–Lévy model $H$, the following are true:
	 	\begin{itemize}
	 		\item $\RR = \bigcup_{a\in [A]^{<\omega}} (\RR \cap L[a])$, 
	 		\item  there is no well-ordering of the reals, and
	 		\item $ A $ has no countable subset.
	 	\end{itemize}
	 \end{Teo}
	 
	 \begin{Lema}[\cite{Beriashvili2018}] \label{lemma: prereq H has the well orders of L[a]}
	 	Consider the model $H$.
	 	Fix an enumeration of the rudimentary functions, and for any $a\in [A]^{<\omega}$ consider the natural order on $a$ as a finite subset of reals. 
	 	Then this fixes a global order $<_a$ in $L[a]$. 
	 	In other words, the relation consisting on triples $(a,x,y)$ such that $x<_a y$ is definable over $H$. 
	 \end{Lema}
	 
	 \begin{Teo}\label{th: prereq paradoxical sets in H}
	 	In $H$, there is
	 	
	 	\begin{enumerate}
	 		\item a Luzin set \cite[Section II]{Pincus1975},
	 		\item \emph{no} Sierpi\'{n}ski set \cite[Theorem 1.6]{Beriashvili2018},
	 		\item a Bernstein set \cite[Theorem 1.7]{Beriashvili2018},
	 		\item  a Vitali set \cite[II.3]{Pincus1975} \footnote{Pincus and Prikry attribute it to Feferman.},
	 		\item a Hamel basis  \cite[Theorem 2.1]{Beriashvili2018}, and
	 		\item a Mazurkiewicz set. \cite[Corollary 0.3]{Beriashvili2022}.
	 	\end{enumerate}
	 	
	 \end{Teo}
	 
	 To construct the PUCs in the models considered, we will need to avoid some geometric obstacles, and this will be done with the help of some results about the transcendence degree of some set theoretical subfields of the reals.
	 In particular, we need to use the following result.

	 	\begin{Lema}[Folklore] \label{lemma: trascendence degree one cohen}
	 	Let $M$ be a model of ZFC and $g$ be a $\mathbb{C}$-generic filter over $M$. 
	 	In $M[g]$, the transcendence degree of $\mathbb{R}$ over the algebraic closure relative to $\mathbb{R}$ of $\mathbb{R} \cap M$ is maximal (i.e. is the cardinality of $\mathbb{R}$). 
	 \end{Lema}

	\section{General setup} \label{section: gral setup}

	To analyze the relationship between paradoxical sets and the Axiom of Choice, we are interested in results of the form ``there is a model of $\mathsf{ZF}$ + there exist  $P$ + no $C$'', where $P$ is some notion of paradoxical set and $C$ is a certain choice principle. 
	When $C$ is the principle of existence of a well-order of the reals, the proofs will have the same structure, which we develop in this section. 
	
	Each model will be an inner model of $V[g][h]$ where $g$ is a $\mathbf{Q}$-generic filter over $V$, and $h$ is a $\mathbf{P}$-generic filter over $V[g]$.
	Usually $\mathbf{P}$ will be a forcing notion approximating the paradoxical set considered, $\mathcal{P}=\cup h$ will be the paradoxical set added by $\mathbf{P}$, and $\mathbf{Q}$ will be an adequate forcing that adds reals, for example, the forcing adding $\aleph_1$-many Cohen reals using finite support.
	
	Theorem \ref{th: general set up}  is built on the work of Brendle et al. (see \cite[Lemma 5.1]{Brendle2018}).
	We want to use the structure of that proof, but write it for a more general set up.
	Theorem \ref{th: general set up} will be the result in full generality, but its corollary, Theorem \ref{th: partition no wo}, will be the result that will be more easily applicable to our purposes. 
	
	\begin{Def}\label{def: real absolute}
		A forcing notion $\mathbf{P}$ is \concept{real absolute} if it is absolute, each condition is a subset of the reals, and the order $\leq_{\mathbf{P}}$ is subset of the order given by $\supseteq$. 
		Namely, there are formulas $\psi$ and $\psi'$ absolute between inner models, such that 
		\begin{align*}
			p\in \mathbf{P} & \iff p\subseteq \RR \text{ and } \psi(p) \\
			\forall p_1, p_2 \in \mathbf{P}, p_1\leq_{\mathbf{P}} p_2  & \iff p_1\supseteq p_2 \text{ and } \psi'(p_1, p_2)
		\end{align*}
	\end{Def}

	\begin{Remark}
		If a forcing $\mathbf{P}$ in $M$ is real absolute and $N$ is an inner model of $M$, then
		\[ \mathbf{P}^{N} = \mathbf{P}^{M} \cap N. \]
	\end{Remark}
	
	\begin{Def}\label{def: real alternating}
		Let $\mathbf{Q}$ be a forcing notion in $V$, and let $\mathbf{P}$ be a forcing notion in $V[g]$ where $g$ is a $\mathbf{Q}$-generic filter over $V$. 
		Then we say that $\mathbf{P}$ and $\mathbf{Q}$ are \concept{real alternating} if the following conditions hold (in $V[g]$): 	
		\begin{enumerate}
			\item for all $p\in \mathbf{P}$, $p\subseteq\RR^{V[g]}$, and there is $r \in \mathbb{R}^{V[g]}$ such that $p\in V[r]$ and $r$ can be computed from finitely many elements of $p$;  and \label{item: item 1 def real alternating}
			\item for all $r\in \mathbb{R}^{V[g]}$, for all $p\in \mathbf{P}$, there is $\bar{p}\in \mathbf{P}$ such that $\bar{p}\leq p$ and $r\in V[\bar{p}]$.
		\end{enumerate} 
	\end{Def}     
	
	\begin{Remark}
		Let $p$ and $r$ be as in item \ref{item: item 1 def real alternating} of Definition \ref{def: real alternating}. 
			Notice that $V[r]$ is a forcing extension of $V$ and a ground of $V[g]$ by Theorem \ref{th: prereq solovay}, since $r\subseteq \omega \subseteq V$ and $r\in V[g]$. 
		Since $r$ can be computed from finitely many elements of $p$ then $r$ will be an element of any model containing $p$. 
		Therefore $V[r]\subseteq V[p]$, and since $p\in V[r]$, then $V[p]= V[r]$.
	\end{Remark}     
	
	\begin{Def}
		Let $V$ be a model of $\mathsf{ZF}$, and let $V[g_0]$, $V[g_1]$, and $V[g_2]$ be three forcing extensions of $V$, not necessarily obtained by the same forcing notion. 
		We say that $\left(V[g_0], V[g_1], V[g_2]\right)$ is a \concept{real bifurcation} if 
		
		\begin{enumerate}
			\item $\mathbb{R}^{V[g_o]}\subsetneq \mathbb{R}^{V[g_1]}$, $\mathbb{R}^{V[g_0]}\subsetneq \mathbb{R}^{V[g_2]}$, and
			\item $V[g_0]=V[g_1]\cap V[g_2]$.
		\end{enumerate}
	\end{Def}
	
\begin{figure}[ht]
		\centering
		\includegraphics[width=1\linewidth]{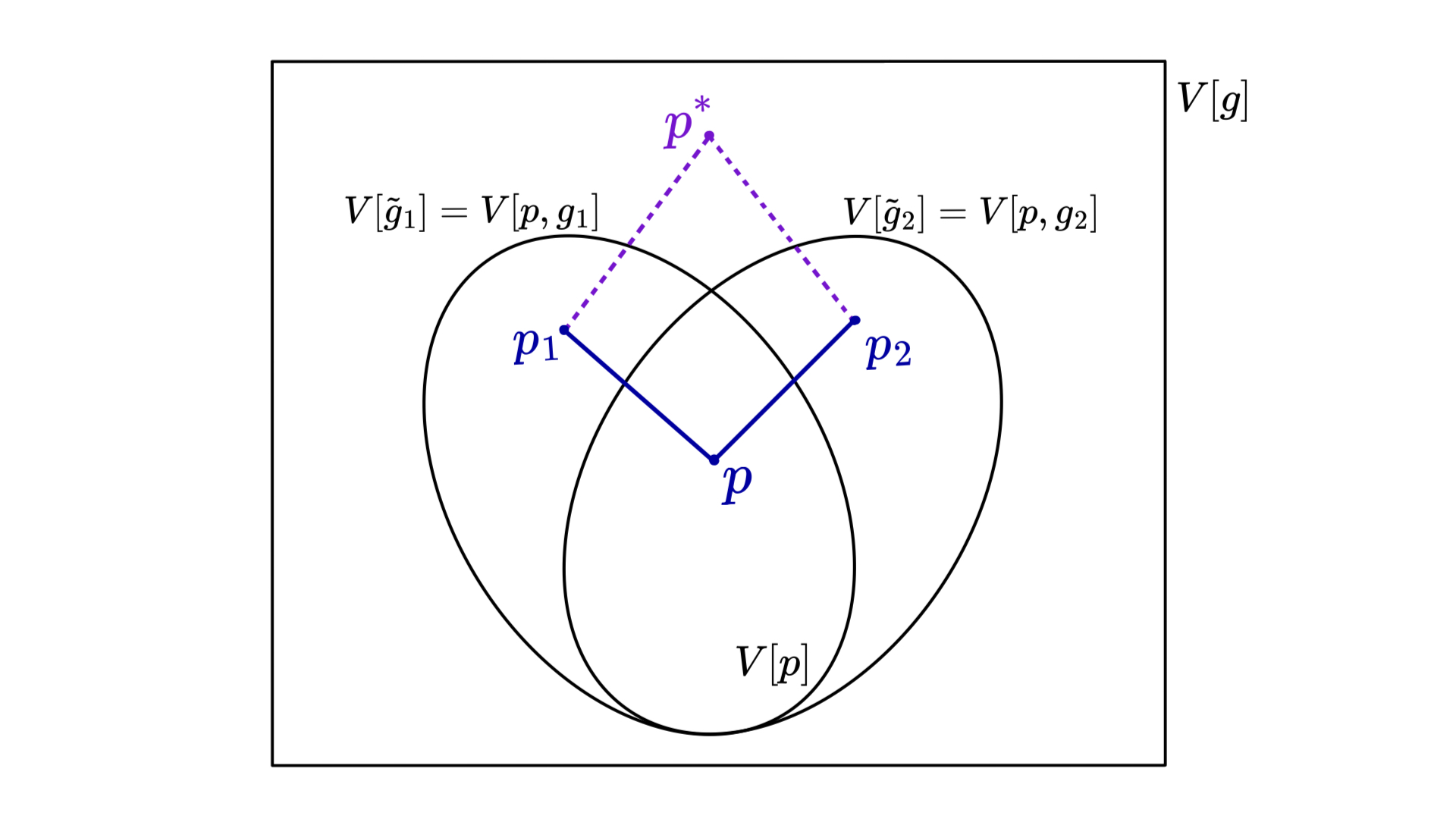}
		\caption{The condition $p\in \mathbf{P}$ is $\mathbf{Q}$-balanced. The compatibility of $p_1$ and $p_2$ is witnessed by $p^\ast$.}
		\label{fig: def-balanced}
	\end{figure}

	\begin{Def}\label{def: balanced}
		Let $\mathbf{Q}$ be a forcing notion in  $V$, let $g$ be a $\mathbf{Q}$-generic filter over $V$ and $\mathbf{P}$ a real absolute forcing notion such that $\mathbf{P}$ and $\mathbf{Q}$ are real alternating.
		We say $\mathbf{P}$ is \concept{$\mathbf{Q}$-balanced over $V$} (in $V[g]$) iff the following statement holds in $V[g]$: 
		\\
		
		For densely many $p\in \mathbf{P}^{V[g]}$, there exist $g_1$, $g_2$ (in $V[g]$) both $\mathbf{Q}$-generic over $V[p]$ such that
		\begin{enumerate}
			\item $V[\tilde{g}_i]=V[p, g_i]$ is a $\mathbf{Q}$-ground for $i=1,2$;
			\item $(V[p], V[\tilde{g}_1], V[\tilde{g}_2])$ is a real bifurcation; and
			\item for all $p_1\in \mathbf{P}^{V[\tilde{g}_1]}$, $p_2 \in \mathbf{P}^{V[\tilde{g}_2]}$ extending $p$, $p_1$ and $p_2$ are compatible. \label{item: balanced- amalgamation}
		\end{enumerate}
		In this case we say that $p$ is a \concept{$\mathbf{Q}$-balanced condition}. 
		See the representation of this situation in Figure \ref{fig: def-balanced}.
	\end{Def} 
	
	Definition \ref{def: balanced} is based on the definition of \emph{balanced} of the book Geometric Set Theory \cite[Definition 5.2.1 and Proposition 5.2.2]{Larson2020} (hence the name). 
	There are differences, the biggest one is that we only require this amalgamation property (item \ref{item: balanced- amalgamation}) for only one pair $g_1$ and $g_2$ instead of \emph{for all}, but the spirit is the same.
	
	\begin{Lema}\label{lemma: Q ground}
	Let $\mathbf{Q}$ be a forcing notion in  $V$, let $g$ be a $\mathbf{Q}$-generic filter over $V$ and $\mathbf{P}$ a real absolute forcing notion such that $\mathbf{P}$ and $\mathbf{Q}$ are real alternating and $\mathbf{P}$ is $\mathbf{Q}$-balanced over $V$.
		
		Let $p\in \mathbf{P}$ be $\mathbf{Q}$-balanced. 
		Suppose $\mathbf{Q}\times \mathbf{Q} \cong \mathbf{Q}$. 
		Then $V[p]$ is a $\mathbf{Q}$-ground of $V[g]$. 
		Moreover, if $\mathbf{Q}$ is homogeneous, for any $g'$ $\mathbf{Q}$-generic filter over $V[p]$, there are densely many conditions $\bar{p}$ in $\mathbf{P}$ such that $\bar{p}$ is $\mathbf{Q}$-balanced (in $V[g']$) and $V[\bar{p}]$ is a $\mathbf{Q}$-ground of $V[g']$. 
	\end{Lema}

	\begin{Dem}
		By definition of balanced, there is a $\mathbf{Q}$-generic filter $g_1$ over $V[p]$. 
		In other words, $V[p]$ is a $\mathbf{Q}$-ground of $V[g_1]$. 
		Also, $V[p,g_1]$ is a $\mathbf{Q}$-ground of $V[g]$. 
		Therefore, $V[p]$ is a $\mathbf{Q}\times \mathbf{Q}$-ground of $V[g]$.
		Since $\mathbf{Q}\times \mathbf{Q} \cong \mathbf{Q}$, we obtain that $V[p]$ is a $ \mathbf{Q}$-ground of $V[g]$.
		\\
		
		For the second part notice that: 
		\[V[g]\models \forall p \in \mathbf{P} \, \exists \bar{p}\leq p \text{ such that }  \bar{p} \text{ is } \mathbf{Q}\text{-balanced and }V[\bar{p}] \text{ is a } \mathbf{Q}-\text{ground}.    \]
		Fix some balanced condition $p$ in $V[g]$. 
		In particular, $p$ is a $\mathbf{Q}$-ground of $V[g]$. 
		Then  there is a condition $q\in \mathbf{Q}$ such that
		\[q \forc{\mathbf{Q} }{V[p]}\forall p' \in \mathbf{P} \, \exists \bar{p}\leq p' \text{ such that }  \bar{p} \text{ is } \check{\mathbf{Q}}\text{-balanced and }V[\bar{p}] \text{ is a } \check{\mathbf{Q}}-\text{ground}. \]
		By homogeneity of $\mathbf{Q}$,
		\[\text{$\mathbbm{1}$}  \forc{\mathbf{Q} }{V[p]}\forall p' \in \mathbf{P} \, \exists \bar{p}\leq p' \text{ such that }  \bar{p} \text{ is } \check{\mathbf{Q}}\text{-balanced and }V[\bar{p}] \text{ is a } \check{\mathbf{Q}}-\text{ground}.  \]
		Hence, for any $g'$ $\mathbf{Q}$-generic filter over $V[p]$, in $V[g']$  there are densely many conditions $\bar{p}\in \mathbf{P}$ such that $\bar{p}$ is $\mathbf{Q}$-balanced and $V[\bar{p}]$ is a  $\mathbf{Q}$-ground of $V[g']$.  
	\end{Dem}
	
	Now we are ready to state the main theorem of this section.\\
	
	\begin{Teo}\label{th: general set up}
		Let $V$ be a model of $\mathsf{ZFC}$. 
		Let $\mathbf{Q}$ be a  forcing notion over $V$, and $g$ be a $\mathbf{Q}$-generic filter over $V$.
		Let $\mathbf{P}$ be a forcing notion over $V[g]$, $h$ be a $\mathbf{P}$-generic filter over $V[g]$, and $\mathcal{P}=\cup h$.
		Suppose the following conditions hold:
		\begin{enumerate}
			
			\item \label{item: no wo- Q homogeneous}$\mathbf{Q}$ is homogeneous and $\mathbf{Q}\times \mathbf{Q} \cong \mathbf{Q}$.
			\item\label{item: no wo- P real absolute} $\mathbf{P}$ is real absolute and $\sigma$-closed.
			\item \label{item: no wo- P Q real alternating}$\mathbf{P}$ and $\mathbf{Q}$ are real-alternating, 
			\item \label{item: no wo- P is Q balanced}	$\mathbf{P}$ is  a $\mathbf{Q}$-balanced forcing over $V$.
		\end{enumerate}
		Then
		\[ L(\RR, \mathcal{P})^{V[g,h]} \models \mathsf{DC} + \neg \mathsf{WO}(\RR). \] 
		
	\end{Teo}
	
	\begin{Remark}
		In Theorem \ref{th: general set up}, notice that the definition of balanced already includes real alternation and $\mathbf{P}$ being real absolute. 
		We include it in the hypotheses of this theorem to facilitate the reading of the proof.
	\end{Remark}
	\\
	
	\begin{Dem}
		First, let us show that $L(\RR, \mathcal{P})^{V[g,h]} \models \mathsf{DC}$.
		Let $R\subseteq A\times A$ be a relation on a nonempty set $A$ such that 
		for every $x\in A$ there  is $y\in A$ with $xRy$. 
		Since $R, A\in L(\RR, \mathcal{P})$, there are some (real and ordinal) parameters that are used together with $\mathcal{P}$ to define $R$ and $A$ by some formulas. 
		Let us call this set of finitely many reals and ordinals by $P$.
		We want to show that there is a sequence $\{x_n\}_{n<\omega}$ such that $x_nRx_{n+1}$ for all $n<\omega$.
		
		Work in $V[g,h]$.
		Fix $x_0\in A \in L(\RR, \mathcal{P})^{V[g,h]}$. 
		Then $x_0$ is definable from $\mathcal{P}$, some (finite) reals, and some ordinals as parameters. 
		Consider $x_0$ as a set definable from $\gamma_0$, $z_0$ and $\mathcal{P}$, where $\gamma_0$ is the least ordinal (that encodes finite ordinals) such that (the decoding of) $\gamma_0$ appears as parameters for the least formula $\phi$ that defines $x_0$. 
		This formula would have some real parameters to define $x_0$. 
		Given $\phi$ and $\gamma_0$, there are some real parameters that work to define $x_0$.
		Choose\footnote{Notice that $V[g,h]$ satisfies $\mathsf{ZFC}$.} such real parameters and encode them by one real $z_0$.
		Notice that $x_0$ is now definable from $z_0\in \RR$ and $\mathcal{P}$.
		For each $n>0$, consider the set $\{y\in A\mid x_nRy\}$. 
		Take $x_{n+1}$ an element of this set that can be defined from the least formula and the least ordinals. 
		Pick some reals that make this formula and these ordinals work for a definition of $x_{n+1}$, encode them in one real $z_{n+1}$.
		Similarly, $x_{n+1}$ is definable from $z_{n+1}$, $\mathcal{P}$ and the set of parameters $P$ (to define $A$ and $R$).
		Take $z= \bigoplus_{n<\omega} z_n$. 
		Then $z$ is a real number that encodes the sequence $\{z_n\}_{n<\omega}$. 
		Since $z\in \RR \in L(\RR, \mathcal{P})^{V[g,h]}$, inside $ L(\RR, \mathcal{P})^{V[g,h]}$ we can decode the sequence $\{z_n\}_{n<\omega}$, which together with $\mathcal{P}$ and $P$, helps us define the sequence $\{x_n\}_{n<\omega}$.
		Then $\{x_n\}_{n<\omega}\in L(\RR, \mathcal{P})^{V[g,h]}$ (it is definable from $\mathcal{P}$ and $P\cup \{z\}$), and $x_nRx_{n+1}$ for all $n<\omega$, as we wanted.
		
		Secondly, we want to show that $L(\RR, \mathcal{P})^{V[g,h]}$ does not have a well-ordering of the reals.
		Suppose the contrary, i.e. that there is some $\mathbf{P}$-generic filter $h$ over $V[g]$, and a formula $\phi(\cdot, \cdot, \vec{x}, \vec{\alpha}, \mathcal{P})$ with $\vec{x}$ a finite sequence in $\mathbb{R}^{V[g,h]}$ and $\vec{\alpha}$ a finite sequence of ordinals such that
		
		\begin{equation}
			V[g,h] \models \phi(\cdot, \cdot, \vec{x}, \vec{\alpha}, \mathcal{P}) \text{ defines a well-ordering of } 2^\omega.
		\end{equation}
		
		Then there is a condition $p_0\in h$ that forces such statement, namely, 
		
		\begin{equation}\label{eq: gral setup step2}
			p_0 \forc{\mathbf{P}}{V[g]} \phi( \cdot, \cdot, \check{\vec{x}}, \check{\vec{\alpha}}, \dot{\mathcal{P}} ) \text{ defines a well-ordering of } 2^\omega,
		\end{equation}
		where $\dot{\mathcal{P}}$ is the name given by $\dot{\mathcal{P}}=\{\tup{p, \check{x}} \mid x\in p \text{ and } p \in \mathbf{P}\}$.
		Notice that we can write $\check{\vec{x}}$ in \ref{eq: gral setup step2} because $\mathbf{P}$ does not add reals, since it is $\sigma$-closed.
		
		On the other hand, since $\mathbf{P}$ is real absolute, there is a formula $\psi$ such that $p\in P \leftrightarrow \psi(p)$ and $\psi$ is absolute between inner models. 
		From now onwards, every time we write $\mathbf{P}$, we are actually interpreting the formula $\psi$ in the corresponding model $V[\cdot]$, which is nothing more than $\mathbf{P}^{V[g]}\cap V[\cdot]$ by absoluteness of $\psi$. 
		Also, when we write $\dot{\mathcal{P}}$, we mean the formula defining  $\dot{\mathcal{P}}=\set{ \tup{p,\check{x}} \mid x\in p \text{ and }\psi(p)}$. 
		
		Now, since $\mathbf{P}$ and $\mathbf{Q}$ are real-alternating, there is $r\in \RR$ such that $p_0\in V[r]$. 
		Let us write $\vec{x}$ as $(x_0, \dots, x_{m-1})$, where $m<\omega$.
		Take $s=r\oplus \bigoplus_{i\in m}  x_i $, $s$ is a real in $V[g]$.
		Again by the property of being real alternating, there is $\bar{p}\leq p_0$ such that $s\in V[\bar{p}]$. 
		Because $\mathbf{P}$ is $\mathbf{Q}$-balanced, there is $p\leq \bar{p}$ which is a balanced condition.
		By real alternation and real absoluteness, we obtain that $V[p]\supseteq V[\bar{p}]$. 
		
		Notice that then $s, r, \vec{x}, p_0 \in V[p]$. 
		From Equation \ref{eq: gral setup step2}, we can write
		\begin{equation}\label{eq: gral setup step3}
			p \forc{\mathbf{P}}{V[g]} \phi( \cdot, \cdot, \check{\vec{x}}, \check{\vec{\alpha}}, \dot{\mathcal{P}} ) \text{ defines a well-ordering of } 2^\omega,
		\end{equation}
		By Lemma \ref{lemma: Q ground}, $V[p]$ is a $\mathbf{Q}$-ground of $V[g]$.
		Let $g'$ be a $\mathbf{Q}$-generic filter over $V[p]$ such that $V[g]=V[p][g']$, and observe that Equation \ref{eq: gral setup step3} is a statement in $V[g]=V[p][g']$. 
		By definability of forcing, we can write this statement as a formula with parameters $\mathbf{P}, p, \check{\vec{x}}, \check{\vec{\alpha}}, \dot{\mathcal{P}}$:
		
		\[V[p, g'] \models \Phi \left(p,\mathbf{P},\check{\vec{x}}, \check{\vec{\alpha}}, \dot{\mathcal{P}} \right)\]
		where $\Phi$ is the formula given by
		
		\[\Phi(\cdot) \iff  p \forc{\mathbf{P}}{V[g]} \phi( \cdot, \cdot, \check{\vec{x}}, \check{\vec{\alpha}}, \dot{\mathcal{P}} ) \text{ defines a well-ordering of } 2^\omega. \]
		
		The statement $\Phi$ has to be forced over $V[p]$ by some condition in $\mathbf{Q}$. 
		Because $\mathbf{Q}$ is homogeneous and all the variables are definable or check names, we get that $\mathbbm{1}_\mathbf{Q}$ already forced it:
		
		\[ \text{$\mathbbm{1}$} _{\mathbf{Q}} \forc{\mathbf{Q}}{V[p]} \Phi \left(\check{p}, \mathbf{P}, \check{\check{\vec{x}}}, \check{\check{\vec{\alpha}}}, \dot{\mathcal{P}} \right). \]
		
		Namely, 
		
		\begin{equation}
			\text{$\mathbbm{1}$} _{\mathbf{Q}}  \forc{\mathbf{Q}}{V[p]} \check{p} \forc{\mathbf{P}}{V[p, \dot{g}]} \phi( \cdot, \cdot, \check{\check{\vec{x}}},\check{\check{\vec{\alpha}}}, \dot{\mathcal{P}} ) \text{ defines a well-ordering of } 2^\omega.
		\end{equation}
		
		Since $p$ is $\mathbf{Q}$-balanced, there are $g_1$, $g_2$ $\mathbf{Q}$-generic filters over  $V[p]$ such that the real bifurcation $(V[p], V[p, g_1], V[p, g_2])$ has the corresponding property of compatibility of conditions. 
		To shorten the notation, we write $V[p, g_i]= V[\tilde{g}_i]$ for $i=1,2$. 
		Notice that $\tilde{g}_i$ does not need to be $\mathbf{Q}$-generic.
		For $i=1, 2$, we get 
		\begin{equation}
			p \forc{\mathbf{P}}{V[\tilde{g}_i]} \phi( \cdot, \cdot, \check{\vec{x}}, \check{\vec{\alpha}}, \dot{\mathcal{P}} ) \text{ defines a well-ordering of } 2^\omega.
		\end{equation}
		
		Take $h_i$ a $\mathbf{P}$-generic filter over $V[\tilde{g}_i]$ such that $p \in h_i$, for $i=1, 2$. 
		Let $\mathcal{P}_i= (\dot{\mathcal{P}})_{h_i}$.
		Then
		\[ V[\tilde{g}_i, h_i] \models \phi(\cdot, \cdot, \vec{x}, \vec{\alpha}, \mathcal{P}_i) \text{ defines a well-ordering of } 2^\omega.\]  
		
		Remember that $(V[p], V[\tilde{g}_1], V[\tilde{g}_2])$ is a real bifurcation and that $\mathbf{P}$ is $\sigma$-closed.
		Notice that, by homogeneity of $\mathbf{Q}$ and $V[p]$ being a $\mathbf{Q}$-ground of both $V[g]$ and $V[\tilde{g}_i]$, we get that
		\[ (\mathbf{P} \text{ is } \sigma\text{-closed})^{V[g]}\implies (\mathbf{P} \text{ is } \sigma\text{-closed})^{V[\tilde{g}_i]} \]
		
		for $i=1, 2$. Therefore, we obtain
		\[ \RR^{V[p]}\subsetneq \RR^{V[\tilde{g}_i]}=\RR^{V[\tilde{g}_i, h_i]}, \] and
		\[\RR^{V[p]}= \RR^{V[\tilde{g}_1]} \cap \RR^{V[\tilde{g}_2]}.\]
		
		Since $V[\tilde{g}_1]$ and $V[\tilde{g}_2]$ have different sets of reals, the respective well orders have to differ at some point. 
		Namely, there is some $\eta \in \OR$ for which the $\eta^{ \text{th}}$-real  given by $\phi$ is different in each model. 
		We then have some digit $n\in \omega$ in which the respective $\eta^{ \text{th}}$ reals differ. 
		Without loss of generality we can write: 
		\[V[\tilde{g}_i,h_i]\models \text{ the } n^{\text{th}} \text{ digit of the } \eta^{\text{th}} \text{ real given by } \phi \text{ is }  i-1.  \ \]
		
		For $i=1,2$, we can find a condition $p_i\leq p $ in $h_i\subseteq \mathbf{P}$ that forces such a statement. 
		Namely, 
		
		\begin{equation}
			p_i \forc{\mathbf{P}}{V[\tilde{g}_i]} \text{ the } \check{n}^{\text{th}} \text{ digit of the } \check{\eta}^{\text{th}} \text{ real given by $\phi$ is }  \check{(i-1 )}.
		\end{equation}
		
		We are exactly in the situation of the definition of $\mathbf{P}$ being $\mathbf{Q}$-balanced over $V$.
		We get then that $p_1$ and $p_2$ are compatible in $V[g]$. 
		
		To obtain a contradiction, we still have to work a bit more. 
		We could be tempted to say that there is a contradiction already, looking at two compatible conditions that force incompatible statements. 
		However, after a closer look, the conditions are forcing incompatible statements over different models. 
		The rest of the proof consists in fixing this obstacle in order to get the desired contradiction.
		
		By Lemma \ref{lemma: Q ground}, for each $i=1, 2 $ there is $\bar{p}_i\in \mathbf{P}^{V[\tilde{g}_i]}$ such that $\bar{p}_i \leq p_i$ and  $V[\bar{p}_i]$ is a $\mathbf{Q}$-ground of $V[\tilde{g}_i]$.  
		Then,
		
		\begin{equation}
			\bar{p}_i \forc{\mathbf{P}}{V[\tilde{g}_i]} \text{ the } \check{n}^{\text{th}} \text{ digit of the } \check{\eta}^{\text{th}} \text{ real given by $\phi$ is }  \check{(i-1 )}.
		\end{equation}
		
		By homogeneity of $\mathbf{Q}$ we have that, for $i=1,2$:
		\begin{equation}
			\text{$\mathbbm{1}$}  \forc{\mathbf{Q}}{V[\bar{p}_i]} \check{\bar{p}}_i \forc{\mathbf{P}}{V[\bar{p}_i, \dot{g}]}\text{ the } \check{\check{n}}^{\text{th}} \text{ digit of the } \check{\check{\eta}}^{\text{th}} \text{ real given by $\phi$ is }  \check{\check{(i-1 )}}. 
		\end{equation}
		
		Notice that $V[\bar{p}_i]\supseteq V[p_i] \supseteq V[p]$ and $\vec{x} \in V[p]$, so the variables of $\phi$ are check names.
		
		Since $V[\tilde{g}_i]$ is a $\mathbf{Q}$-ground of $V[g]$ as well, and $\mathbf{Q}\times \mathbf{Q} \cong \mathbf{Q}$, $V[\bar{p}_i]$ is also a $\mathbf{Q}$-ground of $V[g]$.  This gives us:
		
		\begin{equation}
			\bar{p}_i \forc{\mathbf{P}}{V[g]}\text{ the } \check{n}^{\text{th}} \text{ digit of the } \check{\eta}^{\text{th}} \text{ real is }  \check{(i-1 )}. 
		\end{equation}
		
		More explicitly, 
		\begin{align*}
			\bar{p}_1  &\forc{\mathbf{P}}{V[g]}\text{ the } \check{n}^{\text{th}} \text{ digit of the } \check{\eta}^{\text{th}} \text{ real is }  \check{0}, \text{ and} \\
			\bar{p}_2  &\forc{\mathbf{P}}{V[g]}\text{ the } \check{n}^{\text{th}} \text{ digit of the } \check{\eta}^{\text{th}} \text{ real is }  \check{1}.
		\end{align*}
		
		Let $p^\ast$ be a witness for the compatibility of $\bar{p}_1$ and $\bar{p}_2$ in $\mathbf{P} \cap V[g]$. 
		Then $p^\ast$ forces contradictory statements. 
		Therefore, there is no well-order of the reals in $L(\RR,\mathcal{P})^{V[g,h]}$.
	\end{Dem}
	
	We are interested in applying Theorem \ref{th: general set up} for specific paradoxical sets that are partitions of euclidean spaces in some way. 
	A PUC is clearly a partition but one can think of a Hamel basis as a partition of the reals in the following sense: 
	each real is \emph{covered} by one finite subset of the Hamel basis which spans it
	and this finite subset is unique. 
	So we can think of the reals partitioned into pieces depending on which subset of the Hamel basis spans the real. 

	We wanted  to state Theorem \ref{th: general set up} in full generality for further applications, but for the purpose of this article we will directly apply Theorem \ref{th: partition no wo} instead. 
	To state it we will need some definitions. 
		
		\begin{Def}\label{def: describes real part}
			We say that a formula $\psi(\cdot)$ \concept{describes a real partition} iff it is of the form 
				\[ \psi(p): p\subseteq \RR^n \wedge
			\left(
			\forall s \in [p]^{<\omega} \psi_1(s) \right) \wedge
			\left(
			\forall r \in \RR^m \exists s\in[p]^{<\omega}  \psi_2(r,s) \right),  \]
			where $n,m <\omega$,  $\psi_1$ and $\psi_2$ are absolute between transitive models of set theory, and $\mathsf{ZFC}\vdash \exists p \psi(p)$.
			
		\end{Def}
		
		\begin{Remark}
In applications, this formula $\psi$ will be the description of the paradoxical set considered, the second part will have to do with a condition of \emph{independence} and the third piece of the formula will be related to \emph{covering}.
		\end{Remark}
		
		In particular, the natural formula defining a partition of unit circles does describe a real partition, since $p$ is a PUC iff: 
		\[\psi(p): p\subseteq \mathcal{C} \wedge \text{any two circles are disjoint} \wedge \text{every point is covered}.\]

		\begin{Def}\label{def: forcing adds a real partition}
		Let $V$ be a model of $\mathsf{ZFC}$. 
		Let $\mathbf{Q}$ be a  forcing notion over $V$ and $g$ be a $\mathbf{Q}$-generic filter over $V$.
		Let $\mathbf{P}$ be a forcing notion in $V[g]$ and let $\psi$ be a formula that describes a real partition. 
		We say that $\mathbf{P}$ \concept{adds a real partition according to $\psi$} iff it is of the form 
		\[p \in \mathbf{P} \iff \exists x \in \RR\; V[x]\models \psi(p), \]
		and for any pair $(x,p)$ as before, $x$ can be computed from finitely many elements of $p$. 
		
	Moreover, we request that $\leq_{P}$ is of the form:
		\[p_0\leq_{\mathbf{P}} p_1\iff p_0\supseteq p_1 \wedge \phi(p_0, p_1)\]
		where $\phi$ is absolute between transitive models. 
	\end{Def}

Recall that by Definition \ref{def: describes real part}, $\mathsf{ZFC}\vdash \exists p\psi(p)$. 
In $V[g]$, for all $x\in \RR$ $V[x]\models \mathsf{ZFC}$ and then $V[x]\models \exists p \psi(p)$, so $\mathbf{P}\neq \emptyset$. 
Also, the requirement of the real $x$ being computable from finitely many elements of $p$ implies that $V[p]=V[x]$.
	
	\begin{Def}\label{def: partial condition and extendability}
		Let $\mathbf{P}$ be a forcing notion that adds a real partition as in Definition \ref{def: forcing adds a real partition}. 
		We say that $p\in \mathbf{P}$ is \concept{partial condition} if there is $x\in \RR^{V[g]}$ such that 
		\[V[x] \models p\subseteq \RR^n,
		\forall s \in [p]^{<\omega} \psi_1(s).\] 
		We say that $\mathbf{P}$ satisfies \concept{extendability} if for any partial condition $p$ in $V[x]$ there is a condition $\bar{p}\in \mathbf{P}$ witnessed by $\bar{x}\in \RR$ such that $\bar{p}\supseteq p$ and $x\in V[\bar{x}]$.
		{ \color{black} Moreover, if $p\in \mathbf{P}$ then we request that $\bar{p}$ such that additionally  $\bar{p}\leq_{\mathbf{P}} p$.}
	\end{Def}

		Notice that the reals associated with $p$ and $\bar{p}$ in \ref{def: partial condition and extendability} may differ. 
		They \emph{will} differ for the cases of $\psi$ describing a Mazurkiewicz set and for a partition into unit circles, but we can take $\bar{p}$ so that $p$ and $\bar{p}$ share the same associated real $x$ for the case of Hamel bases. 
		This is because Hamel bases are exactly the maximal linearly independent sets and every partial condition (i.e. a linearly independent set) is extendable to a Hamel basis. 
	
	\begin{Def}\label{def: amalgamation}
		Let $\mathbf{P}$ be a forcing notion that adds a real partition as in Definition \ref{def: forcing adds a real partition}. 
		We say that $\mathbf{P}$ satisfies \concept{amalgamation} (in $V[g]$) if for densely many $p\in \mathbf{P}$, 
		for any $g_1$, $g_2$ mutually $\mathbf{Q}$-generic over $V[p]$  
		and for all $p_1\in \mathbf{P}\cap{V[p, g_1]}$, $p_2 \in \mathbf{P}\cap{V[p,g_2]}$ such that 
		$p_1\leq_{P} p$ and $p_2\leq_{P} p$, 
		$p_1$ and $p_2$ are compatible. 
	\end{Def}
	
	\begin{Teo}[Corollary of Theorem \ref{th: general set up}] \label{th: partition no wo}
		Let $V$ be a model of $\mathsf{ZFC}$. 
		Let $\mathbf{Q}$ be the finite support product of $\omega_1$-many copies of Cohen forcing, and let $g$ be a $\mathbf{Q}$-generic filter over $V$.
		Let $\mathbf{P}$ be a forcing notion over $V[g]$ that adds a real partition, let $h$ be a $\mathbf{P}$-generic filter over $V[g]$, and let $\mathcal{P}=\cup h$.
		
		If $\mathbf{P}$ is $\sigma$-closed and satisfies extendability and amalgamation, then 
		\[ L(\RR, \mathcal{P})^{V[g,h]} \models \mathsf{ZF}+\mathsf{DC}+\neg \mathsf{WO}(\RR) + \psi(\mathcal{P}). \] 
	\end{Teo}
	
	\begin{Dem}
		First, we need to prove that the hypotheses of Theorem \ref{th: general set up} are satisfied.
		\begin{enumerate}
			\item \emph{$\mathbf{Q}$ is homogeneous and $\mathbf{Q}\times \mathbf{Q} \cong \mathbf{Q}$}:
			In this case, 
			$\mathbf{Q}=\mathbf{C}(\omega_1)$ so it satisfies these properties (see Lemma \ref{lemma: prereq C and C(alpha) homogeneous}).
			\item \emph{$\mathbf{P}$ is real absolute and $\sigma$-closed:}
			It is clear that it is real absolute (see Definition \ref{def: real absolute}), by noticing that $L[x]$ does not change through different models containing the same ordinals and $x$, therefore its theory is absolute as well. 
			$\mathbf{P}$ is $\sigma$-closed by hypothesis. 
			\item \emph{$\mathbf{P}$ and $\mathbf{Q}$ are real-alternating:}
			The first condition of real-alternation (see Definition \ref{def: real alternating}) is true because $\mathbf{P}$ adds a real partition. The second is due to $\mathbf{P}$ satisfying extendability. 
			\item \emph{
				$\mathbf{P}$ is  a $\mathbf{Q}$-balanced forcing over $V$:} 
			By amalgamation, there are densely many $p\in\mathbf{P}$ that have the amalgamation property. 
			For such a $p$, there is a real $x$ such that $V[x]=V[p]$. 
			Because of Lemma \ref{lemma: prereq any real is in an inital segment of C(omega1)}, there is some $\alpha<\omega_1$ such that $x\in V[g\restriction\alpha]$. 
			
			By Lemma \ref{lemma: prereq any real in Q extension is Q ground},  $V[x]$ is a $\mathbf{Q}$-ground of $V[g]$.
			Since $\mathbf{Q}\cong \mathbf{Q}\times \mathbf{Q}$, there are $g_1, g_2$ mutually $\mathbf{Q}$-generic over $V[x]$ such that $V[x, g_1, g_2] =V[g]$. 
			Clearly, the tuple $(V[p], V[p,g_1], V[p,g_2])$ is a real bifurcation. 
			Now, take any $p_1\in \mathbf{P}\cap V[p,g_1]$ and $p_2\in \mathbf{P}\cap V[p,g_2]$. 
			By amalgamation, the conditions $p_1$ and $p_2$ are compatible, and $p$ is a $\mathbf{Q}$-balanced condition (see Definition \ref{def: balanced}).
		\end{enumerate}
		
		Applying Theorem \ref{th: general set up}, we get that
		\[ L(\RR, \mathcal{P})^{V[g,h]} \models \mathsf{ZF+DC} + \neg \mathsf{WO}(\RR). \] 
		
		It is left to prove that \[L(\RR, \mathcal{P})^{V[g,h]} \models \psi(\mathcal{P}).\] 
		For this purpose, work inside $V[g,h]$.
		First, since $h\subseteq \mathbf{P}$, we have that for all $p\in h$ there some $x\in \RR$ such that 
		\[ V[x] \models p\subseteq \RR^n,
		\forall s \in [p]^{<\omega} \psi_1(s)  \wedge
		\forall r \in \RR^m \exists s\in[p]^{<\omega}  \psi_2(r,s).  \]
		
		From this and the fact that being a real number is absolute, we have that $\mathcal{P}=uh\subseteq \RR^n$. 
		Let $s\in \mathcal{P}^{<\omega}$, and let $l$ be such that $s=\{s_0, \dots, s_{l-1}\}$. 
		Then there is a finite set of conditions $p_0, \dots, p_{l-1}$ in $h$ such that $s_i\in p_i$ for all $i\in l$. 
		Since $h$ is a filter, there is $p\in h$ such that $p\leq_{\mathbf{P}} p_i$ for all $i\in l$. 
		In particular, $p_i \subseteq p$ for all $i\in l$, and $s\subseteq p$. 
		Since $p\in \mathbf{P}$, there is some real $x$ such that 
		\[V[x]\models \forall \tilde{s} \in [p]^{<\omega} \psi_1(\tilde{s}) \]
		Notice that $s\in [p]^{<\omega}$ so $\psi_1(s)$ holds in $V[x]$.
		Because $\psi_1$ is absolute between inner models, we have that $\psi_1(s)$ holds in $V[g,h]$ as well as in $L(\RR, \mathcal{P})^{V[g,h]}$.
		
		Secondly, notice that $\RR \cap L(\RR, h)^{V[g,h]} = \RR \cap V[g,h] = \RR \cap V[g]$ since $\mathbf{P}$ is $\sigma$-closed. 
		Fix $r\in \RR ^m$.   
		We claim that the set 
		\[D= \{ p\in \mathbf{P} \mid r\in V[p] \} \]
		is dense. 
		Fix $p\in \mathbf{P}$. 
		There is some real $x$ that witnesses $p$ is a condition. 
		By absoluteness, $p$ is a partial condition in $V[x\oplus r]$. 
		By extendability, there is $\bar{p} \supseteq p$ such that $x\oplus r \in V[\bar{p}]$, which implies $r\in V[\bar{p}]$.
		Since $p$ is a condition, we can assume $\bar{p}\leq_{P} p$.
		Since $r\in V[\bar{p}]$, we have shown that $D$ is dense.
		
		Since $h$ is a generic filter, $D\cap h \neq \emptyset$. 
		Let $p \in D\cap h$. 
		By definition of $\mathbf{P}$ again, we have that 
		\[V[p] \models \forall \tilde{r} \in \RR^m \exists s\in[p]^{<\omega}  \psi_2(\tilde{r},s).\]
		Since $r\in V[p]$ by definition of $D$, 
		\[V[p] \models \exists s\in[p]^{<\omega}  \psi_2(r,s).\]
		By absoluteness, there is $s\in[p]^{<\omega}\subseteq [\mathcal{P}]^{<\omega}$ such that $\psi_2(r,s)$ holds in $V[g,h]$ and also in $L(\RR, \mathcal{P})^{V[g,h]}$. 
		Putting everything together, we have that 
		\[ L(\RR, \mathcal{P})^{V[g,h]} \models \psi (\mathcal{P}),\]
		as we wanted to show.
	\end{Dem}
	
	\begin{Lema}\label{lemma: sigma close for partitions}
		Let $\mathbf{Q}$ and $\mathbf{P}$ be as in Theorem \ref{th: partition no wo}. 
		If $\leq_{\mathbf{P}}=\supseteq \restriction (\mathbf{P}\times\mathbf{P})$, then $\mathbf{P}$ is $\sigma$-closed in $V[g].$
	\end{Lema}
	
	\begin{Dem}
		Let $\{p_n\}_{n<\omega}$ be a sequence of $\leq_{\mathbf{P}}$-decreasing conditions. 
		Let $\{x_n\}_{n<\omega}$ be a sequence of reals such that $V[x_n]\models \psi(p_n)$ for all $n<\omega$. 
		We can do this since $V[g]\models \mathsf{AC}$. 
		Take $x=\bigoplus_{n<\omega} x_n$ and $p=\bigcup_{n < \omega} p_n$. 
		
		We claim $p\in V[x]$ is a partial condition.  
		Clearly, $p\subseteq \RR^n$. 
		Fix $s\in [p]^{<\omega}$. 
		We know that $s$ is finite, $\{p_n\}_{n<\omega}$ is a $\leq_{\mathbf{P}}$-decreasing sequence, and $\leq_{\mathbf{P}}$ is a subset of the reverse inclusion in $\mathbf{P}$. 
		Therefore, there is some $n<\omega$ such that $s\subseteq p_n$. 
		Then $V[x_n]\models \psi_1(s)$ and by absoluteness we get that $V[x]\models \psi_1(s)$. 
		
		By extendability, there is a condition $\bar{p}\in \mathbf{P}$ such that $\bar{p}\supseteq p$. 
		Therefore $\bar{p} \leq_{\mathbf{P}} p_n$ for all $n< \omega$.
	\end{Dem}
		
		\begin{Remark}
			Theorem \ref{th: general set up} produces models without a well-ordering on the reals. 
			If $\psi$ is the definition of a paradoxical set, then
			Theorem \ref{th: partition no wo} not only produces a choiceless model but one in which that paradoxical set exists and it is added by the forcing $\mathbf{P}$.
		\end{Remark}

\section[Partitions of the 3-dimensional euclidean space into unit circles]{Main application: Partitions of $\mathbb{R}^3$ into unit circles }\label{section: PUC}

In  this section we will consider an example of a paradoxical set in order to apply Theorem \ref{th: partition no wo}: a partition of $\RR^3$ into unit circles. 
In Subsection \ref{subsection: PUC lit review}, we gave an overview of similar objects constructed with and without choice.
In Sections \ref{section: PUC forcing} and \ref{section: PUC in the Cohen model} we will show models with this paradoxical set but without a well ordering of the reals. The first model will satisfy $\mathsf{DC}$, and in the second model, $\mathsf{AC}_\omega$ does not hold.

\begin{Def}
	Let $\mathcal{C}$ denote the family of circles
of radius one in $\RR^3$.
	We say that $\mathcal{P}\subseteq \mathcal{C}$ is  a \concept{partition of unit circles (PUC)} if $\mathcal{P}$ consists of disjoint circles that cover $\RR^3$, namely, for all $C_1, C_2 \in \mathcal{P}$ we have that $C_1\cap C_2 =\emptyset$, and $\cup \mathcal{P}= \RR^3$.
\end{Def}

As it was discussed in Subsection \ref{subsection: PUC lit review}, Conway and Croft \cite[Appendix]{Conway1964} mentioned for the first time that this object exists using the Axiom of Choice. 
Actually the result they showed is more general and the existence of a partition of $\RR^3$ is only a comment at the end of the appendix of the paper. 
Here we include the proof only for our case, which we took from Jonsson \cite[Lemma 1.7]{Jonsson1998}.\\

\begin{Teo}[$\mathsf{ZFC}$]\label{th: zfc PUC}
	There is a partition of $\RR^3$ into unit circles.
\end{Teo}

\begin{Dem}
	Let $\{x_\alpha\}_{\alpha<\mathfrak{c}}$ be an enumeration of the points in $\RR^3$. 
	We will recursively define $p_\alpha$ for $\alpha<\mathfrak{c}$. 
	
	For $\alpha=0$, set $p_0=\emptyset$. Suppose that $p_\beta$ is defined for all $\beta< \alpha$. 
	If $\alpha$ is a successor ordinal of the form $\beta +1$ and $x_\beta \in \cup p_\beta$, take $p_{\beta+1}=p_\beta$. If $x_\beta \not \in \cup p_\beta$, we will choose a unit circle $C_{\beta}$ such that 
	$x_\beta \in C_{\beta}$ and $C_{\beta} \cap C= \emptyset$ for every $C \in p_{\beta}$. 
	Supposing we can choose such a circle $C_{\beta}$, we define $p_{\beta +1}=p_{\beta}\cup \{C_{\beta}\}$.
	Finally, if $\alpha$ is a limit ordinal, define $p_\alpha= \bigcup_{\beta <\alpha} p_{\beta}$. 
	
	We need to check that the construction is legit, namely, that we can actually choose such a circle $C_\beta$.
	Since we need $C_\beta$ to have radius 1, we only need to choose a center $o_{\beta}$ of the circle and a vector $n_\beta$ normal to the plane in which $C_\beta$ will be contained. 
	If $o_\beta$ and $n_{\beta}$ are fixed, they determine exactly one unit circle. 
	
	First, choose $n_\beta$ such that the plane $\pi_\beta$ determined by $n_\beta$ and the point $x_\beta$ does not contain any of the circles $\{C_\delta\}_{\delta <\beta}$ in $p_\beta$.
	This is possible because there are less than $|\beta|<\mathfrak{c}$ such planes (at most one per circle) and $\mathfrak{c}$ possibilities to choose $n_{\beta}$. 
	
	Second, notice that we need $C_\beta$ to pass through $x_\beta$ and that 
	implies we need $o_\beta$ to be at distance 1 from $x_\beta$. 
	Since we fixed $n_\beta$, the possibilities for $o_\beta$ are contained in the only unit circle $C$ contained in $\pi_\beta$ with center $x_\beta$. 
	For each $\delta <\beta$, $C_\delta \cap \pi_\beta$ consists of at most two points, so there are at most $|\beta|<\mathfrak{c}$ points to avoid. 
	For each of these points $t$, there are at most two options for $o_\beta$ that we have to discard, because they would give rise to two circles $C_\beta$ that would contain $t$ as Figure \ref{fig:puc-zfc} shows.
	But we have $\mathfrak{c}$ choices for $o_\beta$ so we can choose $o_\beta$ so that $C_\beta\cap C_{\delta}=\emptyset$ for all $\delta<\beta$.

	\begin{figure}
		\centering
		\includegraphics[width=0.97\linewidth]{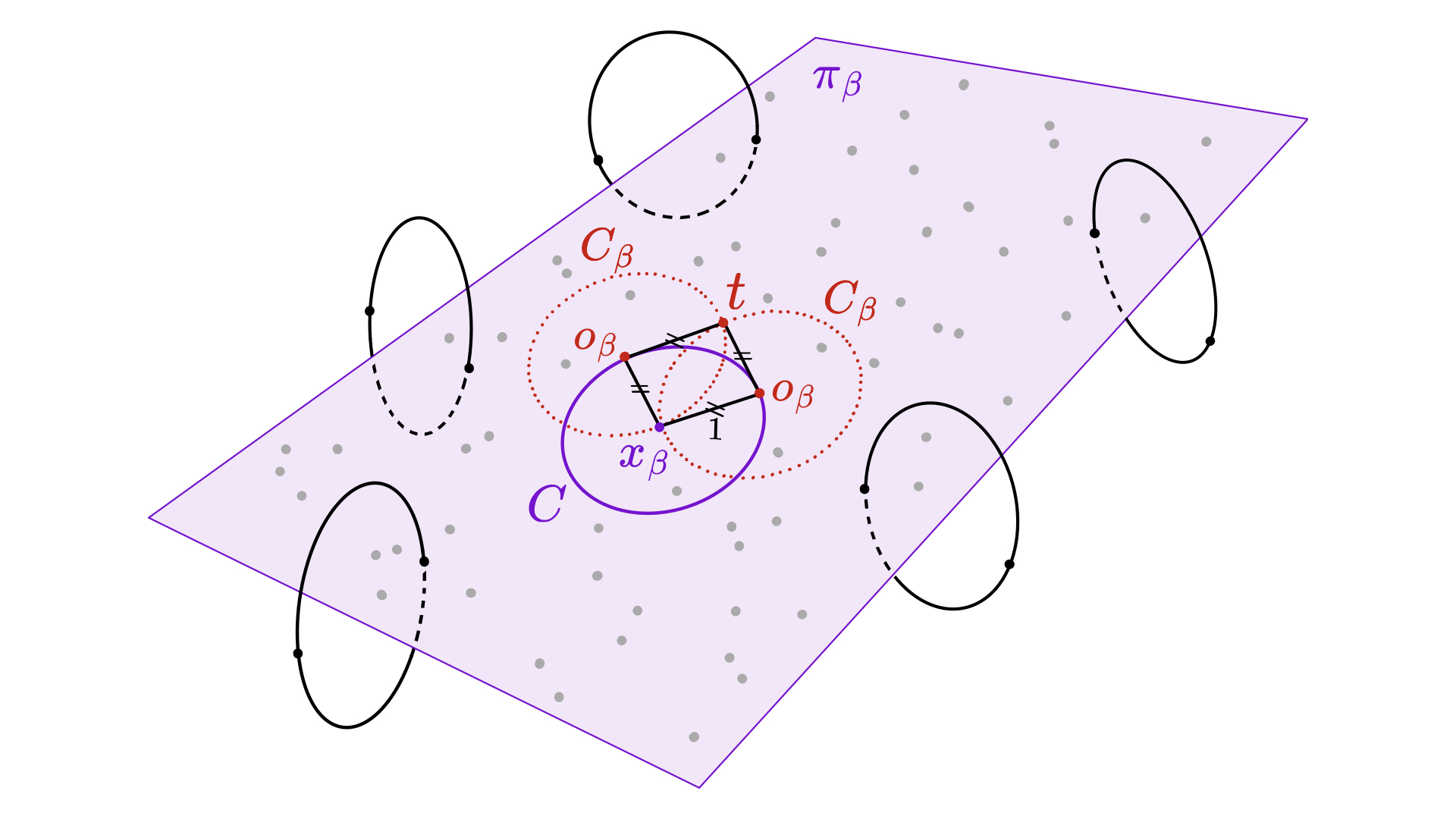}
		\caption{For each point $t$ that we want to avoid, there are two options for $o_\beta$ that we have to discard.}
		\label{fig:puc-zfc}
	\end{figure}

	Take $\mathcal{P}=\bigcup_{\alpha<\mathfrak{c}} p_\alpha$.
	Clearly, $\mathcal{P}$ is a family of unit circles in $\RR^3$.
	For any two circles $C_\beta$ and $C_\alpha$ in $\mathcal{P}$ added in steps $\beta+1$ and $\alpha+1$ of the construction, if $\beta <\alpha$ then $C_\alpha$ was chosen so that $C_\beta \cap C_\alpha =\emptyset$. 
	Moreover, for any $r\in \RR^3$ there is $\alpha < \mathfrak{c}$ such that $r= x_\alpha$. 
	By construction, $r\in \cup p_{\alpha+1} \subseteq \cup \mathcal{P}$. Thus, $\mathcal{P}$ is a partition of $\RR^3$ into unit circles.
\end{Dem}

\subsection{Properties of PUCs}\label{subsection: properties PUCs}

In this subsection we prove properties of PUCs that are known to hold for Mazurkiewicz sets: 
If it is analytic, then it is Borel; and if $V=L$ then there is one which is coanalytic. The latter was independently observed by Linus Richter.

\begin{Prop}\label{prop: PUC if analytic then borel}
	If there is a partition of unit circles that is analytic, it is actually Borel.
\end{Prop}

\begin{Dem}
	Suppose there is a PUC $\mathcal{P}$ which is $\bm{\Sigma^1_1}$.
	We show that its complement is also $\bm{\Sigma^1_1}$ therefore $\mathcal{P}$ is actually Borel. 
	Notice that (unit) circles can be coded by one real.
	\begin{align*}
		C\not\in \mathcal{P} \iff \exists C_1 \exists x_1 \exists C_2 \exists x_2 \text{ such that } & C_1\neq C_2 , C_1, C_2 \in \mathcal{P},\\ 
		& x_1\in C\cap C_1 ; x_2 \in C\cap C_2;
	\end{align*}
	which is again $\bm{\Sigma^1_1}$.
\end{Dem}

\begin{Prop}\label{prop: coanalytic PUC}
If $V=L$, there is a coanalytic partition of $\RR^3$ in unit circles.
\end{Prop}

\begin{Dem}
This is a direct application of \cite[Theorem 1.3]{Vidnyanszky2012}.
\end{Dem}

\subsection{Forcing a PUC} \label{section: PUC forcing}
	
	We aim to apply Theorem \ref{th: partition no wo} to partitions of unit circles. 
	The first challenge we have to face is verifying \emph{extendability}. Indeed, it is not true that any family of disjoint circles (partial conditions) is extendable to a partition inside the same model, so extendability is not trivial (compare with the case of Hamel bases). 
	Furthermore, we have to work much more to show \emph{amalgamation}. This property does not hold if we only consider a forcing poset ordered by reverse inclusion, as it will be the case for Mazurkiewicz sets and Hamel bases (see Section \ref{section: appendix}). 
	Nevertheless, we will be able to show that there is a model of $\mathsf{ZF}+\mathsf{DC}$ with no well order of the reals in which there is a partition of $\RR^3$ into unit circles in Theorem \ref{th: PUC - no wo}.
	We will start by setting some notation and defining the forcing that we will need to construct such a model.

	\begin{Not}
		Given $r\in \RR^n$ we write $\coor(r)$ for the (unordered) set of coordinates. 
		We will extend this notation for circles and planes.
		
		We will think of a circle $C$ as given by parameters $(o, n)$, where $o\in\RR^3$ is its center  and $n \in \RR^3$ is a normal vector of the unique plane that contains $C$. 
		If we choose the normal vectors to be inside the set 
		\[S=\{ (x,y,z)\in \RR^3 \mid x^2+y^2+z^2=1 \text{ and } (z>0 \lor (z=0 \, \land y>0)) \} \cup \{(0,0,1)\}, \]
		then the assignment of a normal vector to any given plane is \emph{unique}. 
		Therefore any circle $C$ has exactly one representation by parameters $(o,n)\in \RR^3\times S$. 
		
		If $C$ is a circle with parameters $(o,n)\in \RR^3 \times S$, then we write $\coor(C)$ for $\coor(o)\cup \coor(n)$.
		Given a model $M$,  we will use ``$C\in M$'' as shorthand for ``$\coor(C)\in M$''.
		
		Consider the set
		\[\overline{S}= \{(a,b,c,d)\in \RR^4 \mid a=1 \lor ( a=0 \land (b=1 \lor (b=0 \land c=1)))\}\]
		Then every plane $\pi$ can be represented \emph{uniquely} by parameters $(a,b,c,d)\in \overline{S}$ such that 
		\[\pi= \{(x,y,z)\in \RR^3 \mid ax+by+cz+d=0\}.  \]
		In this case, we write $\coor(\pi)=\{a,b,c,d\}$. 
		Similarly, ``$\pi \in M$'' is shorthand for ``$\coor(\pi)\in M$''.
		
		If $R$ is a set of circles, planes, or points, we write $\coor(R)$ to denote $\bigcup\{\coor(r)\mid r\in R\}$. 
	\end{Not}

	\begin{Def}\label{def: forcing PUC}
		Let $V$ be a model of $\mathsf{ZFC}$.
		Let $\mathbf{Q}$ be the finite support product of $\omega_1$-many copies of Cohen forcing. 
		Let $g$ be a $\mathbf{Q}$-generic filter over $V$. In $V[g]$, we define a partial order \(\mathbf{P_C}\) as follows:
		\begin{itemize}
			\item \(p\in \mathbf{P_C}\) iff \(\exists x\in \mathbb{R}\) such that  \(V[x]\models p\) is a PUC.
			\item \(p\leq_{\mathbf{P_C}}q\) iff 
			\begin{enumerate}[i.]
				\item $q\supseteq p$, 
				\item there are reals $x$ and $y$ such that  $p$ is a PUC in $V[x]$, $q$ is a PUC in $V[y]$, and $x\in V[y]$. \label{item: ii of def forcing PUC}
				\item $q$ extends $p$ in an \emph{algebraically independent way}. 
				Namely, in $V[y]$, 
				for all $C\in q\backslash p$ with center $o$ and contained in the plane $\pi$, we have that $\pi \not\in V[x]$ and $o\not\in \overline{\RR^{V[x]}(\coor(\pi))}$.
			\end{enumerate}
		\end{itemize}
	\end{Def}
	
			\begin{Not}
		Here, $\overline{F}$ denotes the algebraic closure of $F$ relative to $\RR$ and $\RR^{V[x]}(\coor(\pi))$  is the minimal field containing $\RR^{V[x]}$ and $\coor(\pi)$.
	\end{Not}

	Notice that since ``$\pi \in V[x]$'' means ``$\coor(\pi)\in V[x]$'', and $\coor(\pi)$ is a finite set of reals,  ``$\pi \not \in V[x]$'' means there is at least one coordinate of $\pi$ that is not in $V[x]$. \\
	
 Condition \ref{item: ii of def forcing PUC} in Definition \ref{def: forcing PUC} does not depend on the reals $x$ and $y$. 
	Notice that we can \emph{recover} the real $x$ from the condition $p\in \mathbf{P}$, namely, 
	if $x$ and $x'$ are reals such that $V[x]\models ``p $ is a PUC'' and $V[x']\models ``p $ is a PUC'', then we have that $V[x]=V[x']$. 
	This is due to the fact that the set 
	\[ R= \{ r\in \RR \mid (0,0,r) \in C \text{ where } C \text{ is a circle given by an element of } p  \} \]
	is absolute between models that contain $p$. 
	In each model, $p$ is a PUC and thus covers the respective $z$-axis. 
	Therefore, 
	\[\RR^{V[x]} = R^{V[x]}= R^{V[x']} = \RR^{V[x']}. \]
	Since $x$ and $x'$ are reals, we get $V[x]=V[x']$.\\

	Moreover, $\leq_{P}$ is a partial order. 
	Reflexivity and antisymmetry are clear.
	For transitivity, suppose you have conditions $p,q,r$ witnessed by the reals $x,y,z$ such that $p\leq_{\mathbf{P_C}}q$ and $q\leq_{\mathbf{P_C}}r$. 
	By definition of $p\leq_{\mathbf{P_C}}q$, 
	\[V[y]\models \forall C \in q\backslash p \text{ given by } (o, \pi), \pi \not \in V[x], \text{ and } o\not \in \overline{\RR^{V[x]}(\coor(\pi))}.\]
	By absoluteness, this also is true in $V[z]$.
	Notice that $r\backslash p=r \backslash q \cup q\backslash p$. 
	By definition of $q\leq_{\mathbf{P_C}} r$, 
	\[V[z]\models \forall C \in r\backslash p \text{ given by } (o, \pi), \pi \not \in V[y], \text{ and } o\not \in \overline{\RR^{V[y]}(\coor(\pi))}.\]
	Since $x\in V[y] $, we have $\RR^{V[x]}\subseteq \RR^{V[y]}$, and hence we get that $\pi\not \in V[x]$ and $o\not\in \overline{\RR^{V[x]}(\coor(\pi))}$. 
	\\

	\begin{Remark}
		If $C$ is a unit circle in $V[x]$, its parameters $(o, n)$ are elements of $\RR\cap V[x]$.
		Let $C'$ be the unit circle in $V[y]$ given by $(o,n)$ where $y$ is such that $\RR\cap V[x] \subsetneq \RR \cap V[y]$.
		If we look at $C$ and $C'$ as sets (and not as their definitions) we will get that $C\subsetneq C'$. 
		In other words, the same parameters produce different sets in different models. 
		
		We will alternate between considering the parameters of each circle and the circle itself (the geometrical object) whenever needed, hoping that the reader can perceive whenever this distinction is important. 
	\end{Remark}
	\\
	
	Notice that we can construe $\mathbf{P_C}$ so that $\mathbf{P_C}$ adds a real partition (recall Definition \ref{def: forcing adds a real partition}). 
	We can see the family $\mathcal{C}$ of circles of radii one in $\RR^3$ as the set 
	\[ \mathcal{C}= \{(o,n)\mid o \in \RR^3 \text{ and } n\in S\}=\RR^3\times S.\] 
	Fixing this codification, each condition in $\mathbf{P_C}$ is a subset of $\RR^6$. 
	
	Notice that 
	\[ p \text{ is a PUC } \iff p\subseteq \RR^6 , \forall s \in [p]^2 \psi_1(s) \wedge \forall r\in\RR^3  \exists s \in [p]^{1} \psi_2(r,s), \]
	where 
	\[\psi_1(s) \text{ iff ``} s\subseteq \mathcal{C} \text{ and if }s=\{s_0, s_1\}, \text{ then the circles given by } s_0 \text{ and } s_1 \text{ do not intersect'',}\]
	and
	\[ \psi_2(r,s) \text{ iff ``} s=\{s_0\},  s_0\in \mathcal{C}, \text{ and } r \text{ is covered by the circle given by } s_0\text{''.}\]
	
	First,
	the circles (given by) $s_0=(o_0, n_0)$ and $s_1=(o_1,n_1)$ intersect if and only if the following holds:
	\[ \exists x\in \RR^3 \tupp{o_0-x,n_0}=\tupp{o_1-x,n_1}=0 \text{ and }
	d(x,o_0)=d(x,o_1)=1 \]
	where $\tupp{\cdot,\cdot}$ here denotes the inner product and ``$-$'' is the subtraction of vectors in $\RR^3$. 
	This is a $\bm{\Sigma_1^1}$ property, with parameters $\coor(s_0)\cup\coor(s_1)$.
	By Mostowski’s Absoluteness, it is absolute between transitive models containing the parameters. 

Second, $\psi_2$ is clearly $\bm{\Delta_0}$. 
Furthermore, if $p_1 \leq_{\mathbf{P_C}} p_2$ is given by $p_1\supseteq p_2$ and $\phi(p_1,p_2)$ as in Definition \ref{def: forcing PUC}, then the corresponding $\phi$ is absolute. 

Finally, for every pair $(x,p)$ such that $V[x]\models p \text{ is a PUC}$, there is a circle $C\in p$ which is the only circle in $p$ intersecting the point $(0,0,x)$ and we can compute $x$ from (the parameters of) $C$.
So $\mathbf{P_C}$ adds a real partition.\\

We will construct our model using Theorem \ref{th: partition no wo}
and we will show that $\mathbf{P_{C}}$ satisfies the hypotheses of that theorem one by one, as shown in Lemmas \ref{lemma: PUC extendability}, \ref{lemma: PUC forcing sigma closed} and \ref{lemma: PUC amalgamation}.
Notice that the partial conditions of $\mathbf{P_C}$ (see Definition \ref{def: partial condition and extendability}) are the subsets of $\mathcal{C}\cap V[x]$ for some $x\in \RR^{V[g]}$ which consist of pairwise disjoint circles.
\\

\begin{Lema}\label{lemma: PUC extendability}
	Let $\mathbf{Q}$ be the finite support product of $\omega_1$-many copies of Cohen forcing, let $g$ be a $\mathbf{Q}$-generic filter over $V$. 
	Then $\mathbf{P_{C}}$ in $V[g]$ satisfies extendability. 
\end{Lema}

\begin{Dem}
	Let $p$ be a family of unit circles in $V[x]$ that are pairwise disjoint, where $x\in \RR^{V[g]}$.
	We need to show that we can extend $p$ to $\bar{p}\in V[\bar{x}]$ such that $\bar{p}\supseteq p$ and $x\in V[\bar{x}]$. 
	
Let $\gamma < \omega_1$ be such that $x\in V[g\restriction\gamma]$ (see Lemma \ref{lemma: prereq any real is in an inital segment of C(omega1)}), let $y=\cup (g\restriction\{\gamma\})$, and let $\bar{x}=x\oplus y$. 
	Then $V[\bar{x}]=V[x,y]$, and $y$ is $\mathbf{C}$-generic over $V[x]$.
	We will prove that there is a condition $\bar{p}\in V[\bar{x}]$ such that $\bar{p}\supseteq p$, by strengthening  the construction of a PUC in $\mathsf{ZFC}$ shown in the proof of Theorem \ref{th: zfc PUC}. 
	
	Work in $V[\bar{x}]$. 
	Let $\{x_\alpha\}_{\alpha<\mathfrak{c}}$ be an enumeration of the points in $\RR^3\backslash \cup p$. 
	Here $\cup p$ is the union of all the circles given by $p$ as computed in $V[\bar{x}]$.
	We will recursively define $p_\alpha$ for $\alpha<\mathfrak{c}$. 
	For $\alpha=0$, set $p_0=p$. 
	Notice that $p$ is still a family of disjoint unit circles in $V[\bar{x}]$ by the absoluteness of $\psi_1$. 
	Suppose that $p_\beta$ is defined for all $\beta< \alpha$. 
	If $\alpha$ is a successor ordinal of the form $\beta +1$ and $x_\beta \in \cup p_\beta$ (namely, $x_\beta$ is covered by a circle in $p_\beta$), take $p_{\beta+1}=p_\beta$. If $x_\beta \not \in \cup p_\beta$, we will pick a unit circle $C_{\beta}$ such that $x_\beta \in C_{\beta}$ and $C_{\beta} \cap C= \emptyset$ for all $C \in p_{\beta}$. 
	Assuming we can choose such a $C_\beta$, we define $p_{\beta +1}=p_{\beta}\cup \{C_{\beta}\}$.
	Finally, if $\alpha$ is a limit ordinal, define $p_\alpha= \bigcup_{\beta <\alpha} p_{\beta}$. 
	
	We need to check that the construction is possible, namely, that we can choose such a circle $C_\beta$.
	Again, we only need to choose an origin $o_\beta$ and a normal vector $n_\beta$. 
	
	Notice that if $C\in p\in V[x]$, then its parameters $(o,n)$ would be in $V[x]$. 
	First, we want to choose $n_\beta$ different to all the normal vectors of circles in $p_\beta$.
	Since $|\RR \backslash \RR^{V[x]}|=\mathfrak{c}$, we have in principle continuum many options for $n_\beta$ that are different from all the normal vectors of circles in $p$. 
	Since $p_\beta=p\dot{\cup} \tilde{p}_\beta$ and $|\tilde{p}_\beta|\leq |\beta|$, we have to also avoid choosing at most $|\beta|$-many normal vectors (one per circle in $\tilde{p}_\beta$).
	Because $|\beta|<\mathfrak{c}$, we can choose $n_\beta$ with the property needed. 
	This means that the plane $\pi_\beta$ (determined by $n_\beta$ and $x_\beta$) in which $C_\beta$ will be contained is different from all the planes containing circles in $p_\beta$.
	Therefore, $|\pi_\beta \cap C|\leq 2$ for every circle $C\in p_\beta$.
	
	Secondly, notice that $x_\beta \in C_\beta$ implies that we need $o_\beta$ to be at distance 1 from $x_\beta$. 
	Since we fixed $n_\beta$, the possibilities for $o_\beta$ are contained in the only unit circle $C$ contained in $\pi_\beta$ with center $x_\beta$. 
	Let us choose $o_\beta\in \pi_\beta$ such that at least one coordinate of $o_\beta$ is not in $\overline{F_\beta}$, where 
	\[ F_\beta= \text{ the minimal field containing }(\mathbb{R}\cap V[x]) \cup
	\coor(\tilde{p}_\beta) \cup
	\coor(\pi_{\beta}, x_\beta),\]
	and recall that $\coor(\tilde{p}_\beta)=\bigcup_{\delta<\beta}\coor(o_\delta, \pi_{\delta})$.
	
	We can choose such an $o_{\beta}$ because of the Lemma \ref{lemma: trascendence degree one cohen}, and because we still have one degree of freedom for a point in $ \RR^3$ after prescribing \(o_{\beta} \in \pi_{\beta}\) and \(d (o_{\beta}, x_{\beta})=1\). 
	See Figure \ref{fig:puc-zfc}.
	We might not be able to choose all the coordinates of $o_\beta$ to not be in $\overline{F_\beta}$, but we only need one of them to not be in $\overline{F_\beta}$. 
	
	Let $C_\beta$ be the circle determined by $(o_\beta, n_\beta)$. 
	We need to check that it satisfies the requirements that we requested in the recursive definition. 
	Clearly $x_\beta \in C_\beta$, since \( o_{\beta}, x_{\beta}\in \pi_{\beta}\), \( d(o_{\beta}, x_{\beta})=1 \) and $o_{\beta}$ is the center of $C_{\beta}$. 
	Fix $\tilde{C}\in p_\beta$.
	We want to show $C_\beta \cap \tilde{C}=\emptyset$. 
	Suppose there is some \(t\in C_{\beta}\cap \tilde{C}\).
	If $\tilde{C}\in p$, its parameters belong to \(\mathbb{R} \cap V[x]\).
	We can calculate $o_{\beta}$ from $(o, \pi, x_{\beta}, \pi_{\beta})$ ``algebraically'': $t$ is one of the (at most two) intersection points of the only unit circle $\tilde{C}$ given by $(o, \pi)$ and the plane $\pi_\beta$, and $o_\beta$ is then one of the (at most two) points in $\pi_\beta$ such that $d(o_\beta,x_\beta)=d(o_\beta,t)=1$.
	See Figure \ref{fig: PUC-fig3}.
	Moreover, in such a situation, there are polynomials $P_i$ of degree $4$ with coefficients in the minimal field containing \( \coor(o, \pi, x_{\beta}, \pi_{\beta})\) such that \(P_i(o^{(i)}_{\beta})=0\) for $i=1,2,3$. 
	Here $o_{\beta}^{(i)}$ denotes the $i^{\text{th}}$ coordinate of the point $o_\beta$.
	This implies that all the coordinates of $o_\beta$ belong to $\overline{F_\beta}$, contradicting the choice of $o_{\beta}$.

	\begin{figure}[ht]
		\centering
		\includegraphics[width=0.8\linewidth]{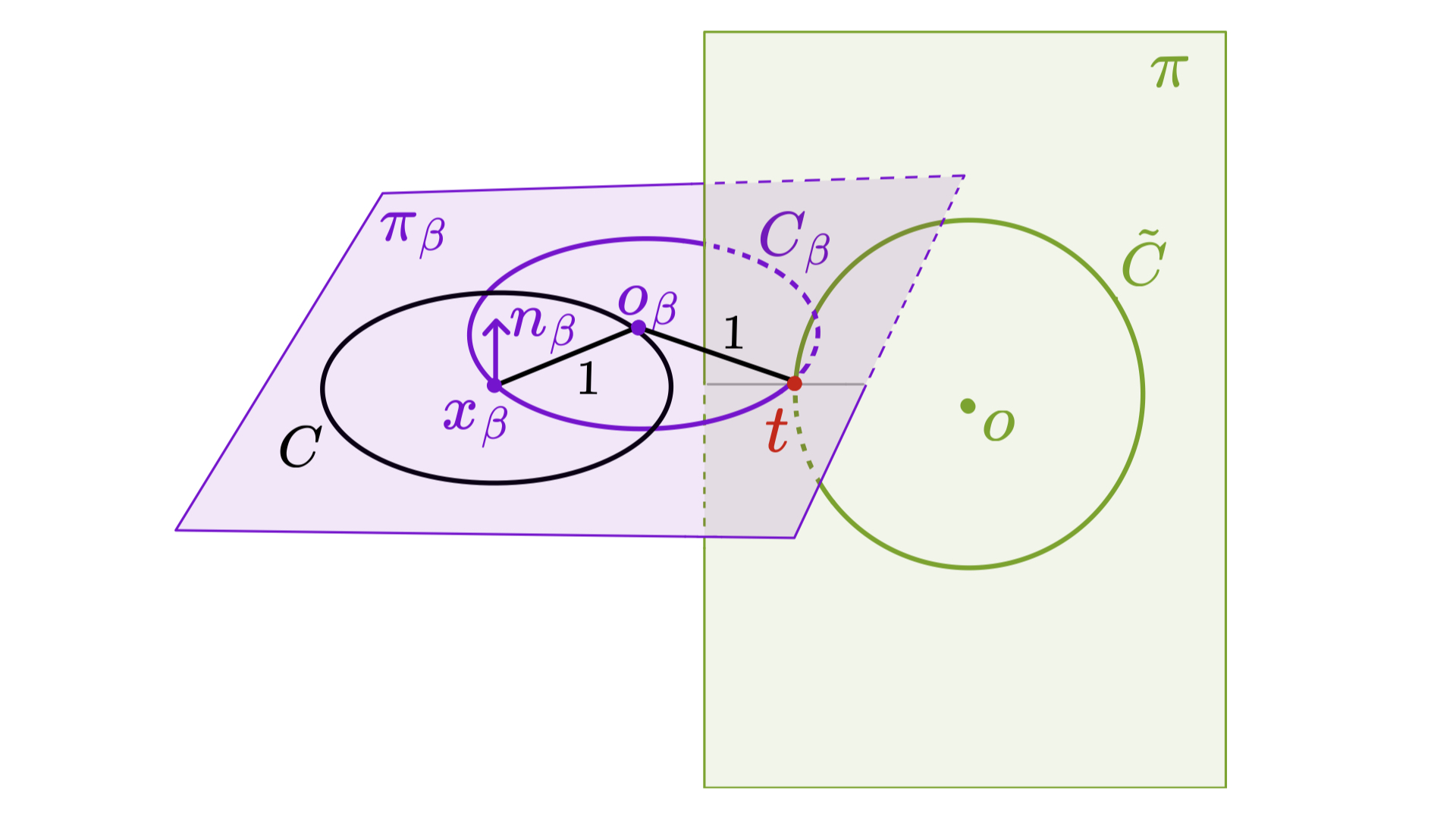}
		\caption{$t$ is one of the (at most two) intersection points of the only unit circle $\tilde{C}$ given by $(o, \pi)$ and the plane $\pi_\beta$, and $o_\beta$ is then one of the (at most two) points in $\pi_\beta$ such that $d(o_\beta,x_\beta)=d(o_\beta,t)=1$.} 
		\label{fig: PUC-fig3}
	\end{figure}

	The case in which $\tilde{C}\in \tilde{p}_\beta$ is analogous. 
	$\tilde{C}$ must have been added in some step $\delta +1$.
	We obtain a contradiction from 
	\(\tilde{P}_i(o^{(i)}_{\beta})=0\) for $i=1,2,3$; where $\tilde{P}_i$ is some polynomial that has coefficients in the minimal field containing
	\( \coor(o_{\delta}, \pi_{\delta}, p_{\beta}, \pi_{\beta}).\)\\

	Take $\bar{p}=\bigcup_{\alpha<\mathfrak{c}} p_\alpha$.
	For any two circles $C, D \in \bar{p}$, we want to show that $C\cap D=\emptyset$.
	This is clear for $C,D$ added in step $0$, namely, $C,D\in p$. 
	If they were added in different steps, for example, $D$ strictly after $C$, there is $\alpha<\omega_1$ such that $D=C_\alpha$ and $C\in p_\alpha$.
	By construction, $C_\alpha\cap C =\emptyset$. 
	Moreover, for any $r\in \RR^3$ either $r\in \cup p$ or there is	
	$\alpha < \mathfrak{c}$ such that $r= x_\alpha$. 
	In the first case, $r\in \cup \bar{p}$ since $p=p_0\subseteq \bar{p}$.
	In the second case, $r\in \cup p_{\alpha+1} \subseteq \cup \bar{p}$ by construction. 
	Thus, $\bar{p}$ is a partition of $\RR^3$ into unit circles that extends $p$.
	\\
	
	Additionally, if $p$ is a condition such that $V[x]\models p$ is a PUC, then by construction we have $\bar{p}\leq_{\mathbf{P_C}} p$. 
\end{Dem}

It is not true that every partial PUC in a model $M$ of $\mathsf{ZFC}$ can be extended to a (complete) PUC inside the model $M$. 
	Notice that the proof of Theorem \ref{th: zfc PUC} shows that any family of disjoint unit circles of cardinality less than $\mathfrak{c}$ can be extended to a partition of $\RR^3$ into unit circles. 
This is not true for families of disjoint unit circles that have cardinality $\mathfrak{c}$ even if there are still $\mathfrak{c}$ points to be covered.
For example, a similar proof of Theorem \ref{th: zfc PUC} shows that we can partition $\RR^3\backslash l$ into unit circles, where $l$ is any line in $\RR^3$. 
This is a family of disjoint unit circles, they cover exactly $\RR^3\backslash l$ so there are $|l|=\mathfrak{c}$ points not covered.
Nevertheless, there is no circle that we can add to this family to cover all $\RR^3$.
However, Lemma \ref{lemma: PUC extendability} says that we can always do extend a family of disjoint unit circles to a PUC when we make \emph{more space} for it, namely, add more reals to the model.

\begin{Remark}
	We will use and abuse the notation $p_0\leq_{\mathbf{P}} p_1$ even when $p_0$ and $p_1$ are partial conditions. 
	If the models where we are considering $p_0$ and $p_1$ are fixed, for example $p_0\in V[y_0]$ and $p_1\in V[y_1]$, we can reuse Definition \ref{def: forcing PUC}. 
	One may not be able to recover the real $x$ from a partial condition,  so $\leq_{\mathbf{P_M}}$ is not a relation between partial conditions. 
	We could define it as a relation between pairs $(y,p)$ where $V[y]\models$ ``$ p$ is a family of disjoint unit circles''.  
	In this case, it will not be a partial order because antisymmetry fails, but the relation is transitive by the same argument that shows $\leq_{\mathbf{P_C}}$ as a relation on $\mathbf{P_C}$ is transitive.
	
	Using this notation, the proof of Lemma \ref{lemma: PUC extendability} gives us that for any partial condition $p$ in a model $V[x]$, there is a condition $\bar{p}\in \mathbf{P}$ and $\bar{x}$ such that $x\in V[\bar{x}]$ and 
	$\bar{p}\leq_{\mathbf{P_C}} p$. 
	We will need this for the proof of Lemma \ref{lemma: PUC forcing sigma closed}.
\end{Remark}

\begin{Lema} \label{lemma: PUC forcing sigma closed}
	Let $\mathbf{Q}$=$\mathbf{C}(\omega_1)$. 
	Let $g$ be a $\mathbf{Q}$-generic filter over a model $V$ of $\mathsf{ZFC}$. 
	Then $\mathbf{P_C}$ is $\sigma$-closed in $V[g]$. 
\end{Lema}

\begin{Dem}
	Work in $V[g]$. 
	Let $\{p_n\}_{n<\omega}$ be a sequence of decreasing conditions. 
	Let $\{x_n\}_{n<\omega}$ be a sequence of reals such that $V[x_n]\models ``p_n$ is a PUC'' for all $n<\omega$. 
	We can do this since $V[g]\models \mathsf{AC}$. 
	Take $x=\bigoplus_{n<\omega} x_n$ and $p=\bigcup_{n < \omega} p_n$. 
	By the proof of Lemma \ref{lemma: sigma close for partitions}, $V[x]\models p$ is a family of disjoint unit circles, so $p$ is a partial condition.
	
	Fix $n<\omega$. 
	Then $p \leq_{\mathbf{P}} p_n$: for all $C\in p\backslash p_n$, there is $m>n$ such that $C\in p_m\backslash p_n$.
	Since $p_m\leq_{\mathbf{P}}p_n$,
	$\pi \not \in V[x_n]$, and $o\not\in \overline{\RR^{V[x_n]}(\coor(\pi))}$ in $V[x_m]$.
	This also holds in $V[x]$ by absoluteness.
	
	Using Lemma \ref{lemma: PUC extendability}, we find $\bar{x}\in \RR$ and $\bar{p}\in \mathbf{P_C}$ such that $V[\bar{x}]\models \bar{p}$ is a PUC. 
	By the Remark above we know that $\bar{p} \leq_{\mathbf{P_C}} p$. 
	Since $p \leq_{\mathbf{P_C}} p_n$ for all $n<\omega$, by transitivity we get that $\bar{p}\leq_{\mathbf{P_C}} p_n$ for all $n<\omega$.
\end{Dem}

Lemma \ref{lemma: PUC forcing sigma closed} implies that $\mathbf{P_C}$ does not add reals.
It is very important to have this property for our purposes because, intuitively, we could have been trying to add a PUC $h$ by partial versions of it, while at the end adding new reals which would not have been considered in the partial approximations.
Therefore the forcing would not ensure that $h$ covers all the points in $\RR^3$ in the extension.
\\

Finally, we are ready to prove the last lemma of this section. 

\begin{Lema} \label{lemma: PUC amalgamation}
	Let $\mathbf{Q}$ be the finite support product of $\omega_1$-many copies of Cohen forcing, let $g$ be a $\mathbf{Q}$-generic filter over $V$.	
	Then $\mathbf{P_C}$ satisfies amalgamation in $V[g]$.
\end{Lema}

\begin{Dem}
	We need to prove that for densely many $p\in \mathbf{P_C}$, for any $g_1, g_2$ mutually $\mathbf{Q}$-generic over $V[p]$, and for all $p\in \mathbf{P_C}\cap V[p,g_1]$, $p_2\in \mathbf{P_C}\cap V[p,g_2]$ such that  $p_1, p_2 \leq_{\mathbf{P_C}} p$, we get that $p_1$ and $p_2$ are compatible.
	
We can assume that there are $x,y,z\in \RR^{V[g]}$ such that 
	\begin{align}
		V[x]\models& p \text{ is a PUC}, \label{eq: amalg puc p in V[x]}\\
		V[x,y] \models& p_1 \text{ is a PUC, and} \\
		V[x,z]\models& p_2 \text{ is a PUC}; \label{eq: amalg puc p2 in V[x,z]}
	\end{align}
	and $y$ and $z$ are mutually generic Cohen reals over $V[x]$. 
	
	Work in $V[x,y,z]$.
	It is clear that \(p_1 \cup p_2\) is a family of unit circles. 
	We only need to prove that they are also disjoint to assert that $p_1\cup p_2$ is a partial condition. 
	We already showed that if two circles are disjoint in $V[x]$ then they are disjoint in $V[x,y,z]$. 
	So let \(C_1 \in p_1 \backslash p\), \(C_2 \in p_2\backslash p\), and suppose \(C_1\cap C_2 \neq \emptyset \). 
	Let \(\pi_i, o_i\) be the plane and the origin, respectively, of $C_i$, for $i=1,2$.
	Notice that $\coor(\pi_1, o_1)\subseteq V[x,y]$ and $\coor(\pi_2, o_2)\subseteq V[x,z]$
	We have two cases, given by the cardinality of $C_1\cap C_2$. We aim to reach a contradiction.
	
	\begin{figure}[ht]
		\centering
		\includegraphics[width=0.7\linewidth]{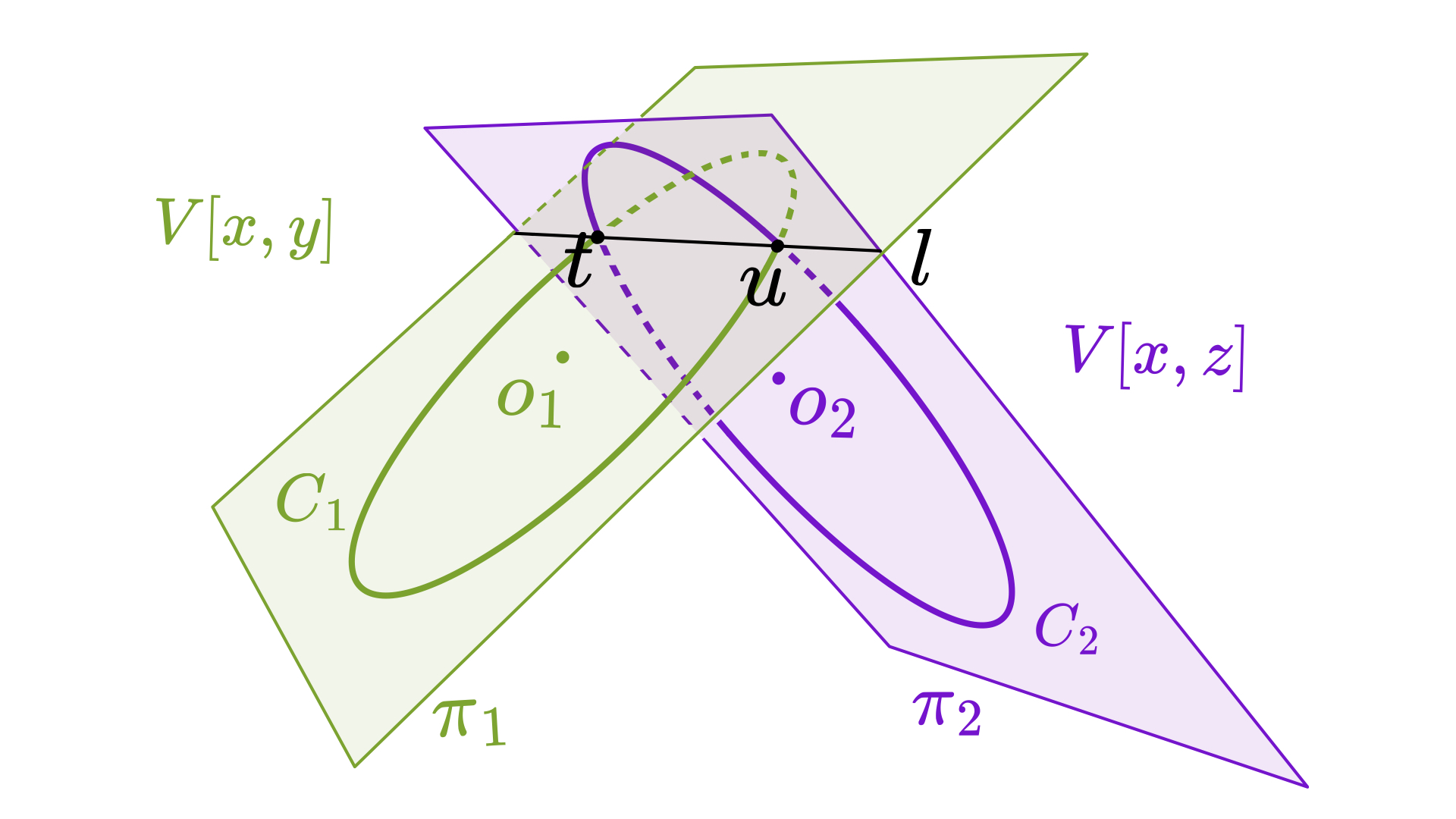}
		\caption{Two circles from different models intersecting in two points $t$ and $u$.} 
		\label{fig: PUC-fig4}
	\end{figure}

	\textbf{Case 1.} Assume \(C_1 \cap C_2 = \{t,u\}\), with $t\neq u$.
	Observe that \(\pi_1 \neq \pi_2\). 
	Otherwise, 
	\[\coor(\pi_1)= \coor(\pi_2) \in V[x,y] \cap V[x,z]= V[x], \]
	which contradicts the requirement $\pi_1 \not \in V[x]$, given by $p_1\leq_{\mathbf{P_C}} p$.
	
	Since \(\pi_1 \neq \pi_2 \), we obtain that  \(C_1 \cap \pi_2 = C_2 \cap \pi_1 =  \{t, u\} \) )(see Figure \ref{fig: PUC-fig4}).
	Now, we can calculate $o_2$ algebraically using $\pi_1, o_1, \pi_2$: 
	we get $t$ and $u$ from computing  \(C_1\cap \pi_2, \) and $C_1$ is given by $\pi_1$ and $o_1$.
	Now, $o_2$ is one of the two points in $\pi_2$ such that \( d(o_2, t)= d(o_2,u) =1\). 
	
	Moreover, in such a situation, there are polynomials $P_i$ of degree $2$ with coefficients in the minimal field containing \(\coor(o_1, \pi_1, \pi_2)\) such that \(P_i(o_2^{(i)})=0\), for $i=1,2,3$. 
	Here $o_2^{(i)}$ denotes the $i^{\text{th}}$ coordinate of the point $o_2$.
	Remember that $p_2\leq_{\mathbf{P_C}} p$ so 
	\( o_2 \not \in \overline{\mathbb{R}^{V[x]}(\coor(\pi_2))}\), 
	namely, there is a coordinate of $o_2$ that does not belong to this field. 
	Suppose without loss of generality that it is $o_2^{(1)}$.
	Take $B\subseteq \{o_2^{(1)}\}\cup \coor(\pi_2) $ maximal such that $B\subseteq \RR^{V[x,z]}$ is algebraically independent over $\RR^{V[x]}$ and contains $o_2^{(1)}$.
	Then, by \cite[Proposition 3.15]{Fatalini2024}, $B$ is also algebraically independent over $\RR\cap V[x,y]$. 
	Recall that $\coor(o_1, \pi_1) \subseteq \RR\cap V[x,y]$. 
	This leads to a contradiction, since \(P_1(o^{(1)}_2)=0 \). 
	
	\textbf{Case 2.} $C_1\cap C_2=\{t\}$.
	
	\textbf{Case 2a.} Suppose that there is a circle $C\neq C_1$ with parameters in $V[x,y]$ such that $C\cap C_2=\{u\}$ and $u\neq t$, 
	as Figure \ref{fig: PUC amalgam two points two circles} shows.
	Let $o$ and $\pi$ be the origin and plane of $C$, respectively.
	Then, similarly to Case 1, we can compute algebraically all the coordinates from $o_2$ using $\coor(o_1, \pi_1, o, \pi, \pi_2)$. 
	The contradiction is analogous. 
	
	\begin{figure}
		\centering
		\includegraphics[width=0.95\linewidth]{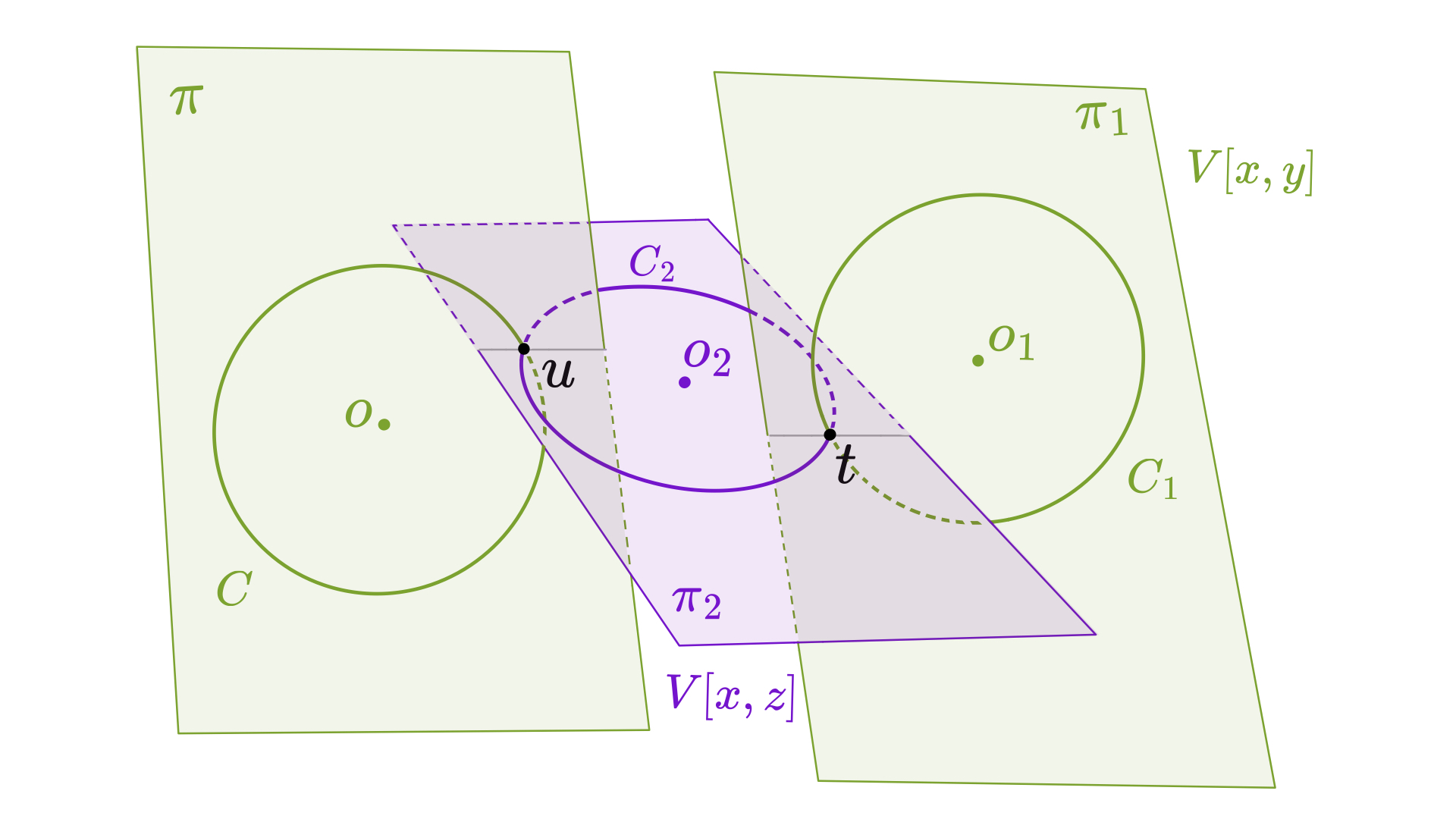}
		\caption{There are circles $C$ and $C_1$ with parameters in $V[x,y]$ such that $C\cap C_2=\{u\}$, $C_1 \cap C_2=\{t\}$, and $u\neq t$.}
		\label{fig: PUC amalgam two points two circles}
	\end{figure}

	\textbf{Case 2b.}  Suppose that there is a circle $C$ with parameters in $V[x,y]$ such that $C\cap C_2=\{t\}$. 
	Then $t\in C \cap C_1$, so $t\in V[x,y]$.
	We know that $t$ is the only point in $C_2\cap V[x,y]$. 
	If not, we could easily define a circle in $V[x,y]$ passing through a possible second point $u$, and this situation was discarded in Case 2a. 
	So we know that
	\[V[x,y,z] \models t \text{ is the only element of } V[x,y] \cap C(o_2,\pi_2), \]
	where $C(o_2, \pi_2)$ describes the unique unit circle with origin $o_2$ contained in the plane $\pi_2$. 
	
	Recall that $y$ is $\mathbf{C}$-generic over $V[x, z]$. 
	There is then a condition $s\in y \subseteq \mathbf{C}$ and $\dot{t}$ a $\mathbf{C}$-name for $t$ in $V[x,z]$ such that 
	\begin{equation} \label{eq: PUC amalgamation 2b}
		s \forc{\mathbf{C}}{V[x,z]} \dot{t}  \text{ is the only element of } V[x,\dot{g}] \cap C(\check{o}_2,\check{\pi}_2),
	\end{equation}
	where $\dot{g}$ is the usual name for the $\mathbf{C}$-generic real $y$. 
	Split $y$ in two mutually generic Cohen reals $y_1, y_2$ according to $s$ as in Definition \ref{def: prereq split of a cohen real}.
	From Equation \ref{eq: PUC amalgamation 2b}, we get that 
	\begin{align*}
		V[x,z,y_1] \models& t_1 \text{ is the only element of } V[x,y_1] \cap C(o_2,\pi_2), \text{ and} \\
		V[x,z,y_2] \models& t_2 \text{ is the only element of } V[x,y_2] \cap C(o_2,\pi_2),
	\end{align*}
	where $t_1= \dot{t}_{y_1}$ and $t_2= \dot{t}_{y_2}$.
	
	Since $V[x,y_1], V[x,y_2]\subseteq V[x,y]$ and $t, t_1, t_2 \in C_2$, we obtain that $t=t_1=t_2$. 
	Then, $t \in V[x,y_1]\cap V[x,y_2]$, so $t\in V[x]$.
	Since $V[x]\models p$ is a PUC, $t$ was covered by some circle $C$ in $p$. 
	Hence, $C\cap C_1 \neq \emptyset$, which contradicts $p_1\leq_{\mathbf{P_C}} p$.
	
	\textbf{Case 2c. } $C_1$ is the only circle (with parameters) in $V[x,y]$ such that $C_1\cap C_2 \neq \emptyset$. 
	
	Similarly to Case 2b, we have that there is an $s\in y$ such that 
	\begin{equation} \label{eq: PUC amalgamation 2c}
		s \forc{\mathbf{C}}{V[x,z]} \tau  \text{ is the only circle from } V[x,\dot{g}] \text{ that intersects } C(\check{o}_2,\check{\pi}_2).
	\end{equation}
	Split $y$ again in two mutually generic Cohen reals $y_1, y_2$ containing $s$. 
	From Equation \ref{eq: PUC amalgamation 2c}, we get that 
	\begin{align*}
		V[x,z,y_1] \models& D_1 \text{ is the only circle from } V[x,y_1] \text{ that intersects } C(o_2, \pi_2), \text{ and} \\
		V[x,z,y_2] \models& D_2 \text{ is the only circle from } V[x,y_2] \text{ that intersects }  C(o_2, \pi_2),
	\end{align*}
	where $D_1= \tau_{y_1}$ and $D_2= \tau_{y_2}$.
	
	Since $V[x,y_1], V[x,y_2]\subseteq V[x,y]$, and $C_1, D_1, D_2 $ define circles that intersect $C_2$, we obtain that $C_1=D_1=D_2$. 
	Then, $C_1 \in V[x,y_1]\cap V[x,y_2]$, so $C_1\in V[x]$, i.e., $\coor(C_1)\in V[x]$.
	Since $C_1\in p_1\backslash p$, by definition of $p_1\leq_{\mathbf{P_C}} p$ we get that $\pi\not \in V[x]$
	This is a contradiction. \\
	
	Taking all the cases into account, we obtain that $p_1\cup p_2$ is a family of disjoint unit circles, and therefore it is a partial condition with respect to $\mathbf{P_C}$.
	Moreover, considering $V[x,y,z]$ as the model containing $p_1\cup p_2$ we claim that $p_1\cup p_2 \leq_{\mathbf{P_C}} p_1, p_2$. 
	If $C\in (p_1\cup p_2) \backslash p_2$, namely, $C\in p_1\backslash p$, we know that $o\not \in \overline{\RR^{V[x]}(\coor(\pi))}$.
	Take $B\subseteq \coor(\pi)\cup \coor(o)$ a maximal algebraically independent set over $\RR^{V[x]}$ containing the coordinate of $o$ that is not in $\overline{\RR^{V[x]}(\coor(\pi))}$. 
	Then $B$ is also algebraically independent over $\RR^{V[x,y]}$ by \cite[Proposition 3.15]{Fatalini2024}, and hence $o\not \in \overline{\RR^{V[x,y]}(\coor(\pi))}$.
	
	Finally, by Lemma \ref{lemma: PUC extendability}, we can obtain a condition $\bar{p}\in V[\bar{x}]$ such that $\bar{p}\leq_{P}p_1\cup p_2$, and such that $V[x,y,z]\subseteq V[\bar{x}]$. 
	By transitivity, $\bar{p} \leq_{\mathbf{P_C}} p_1, p_2$ as we wanted.
\end{Dem}

\begin{figure}[h]
	\centering
	\includegraphics[width=\linewidth]{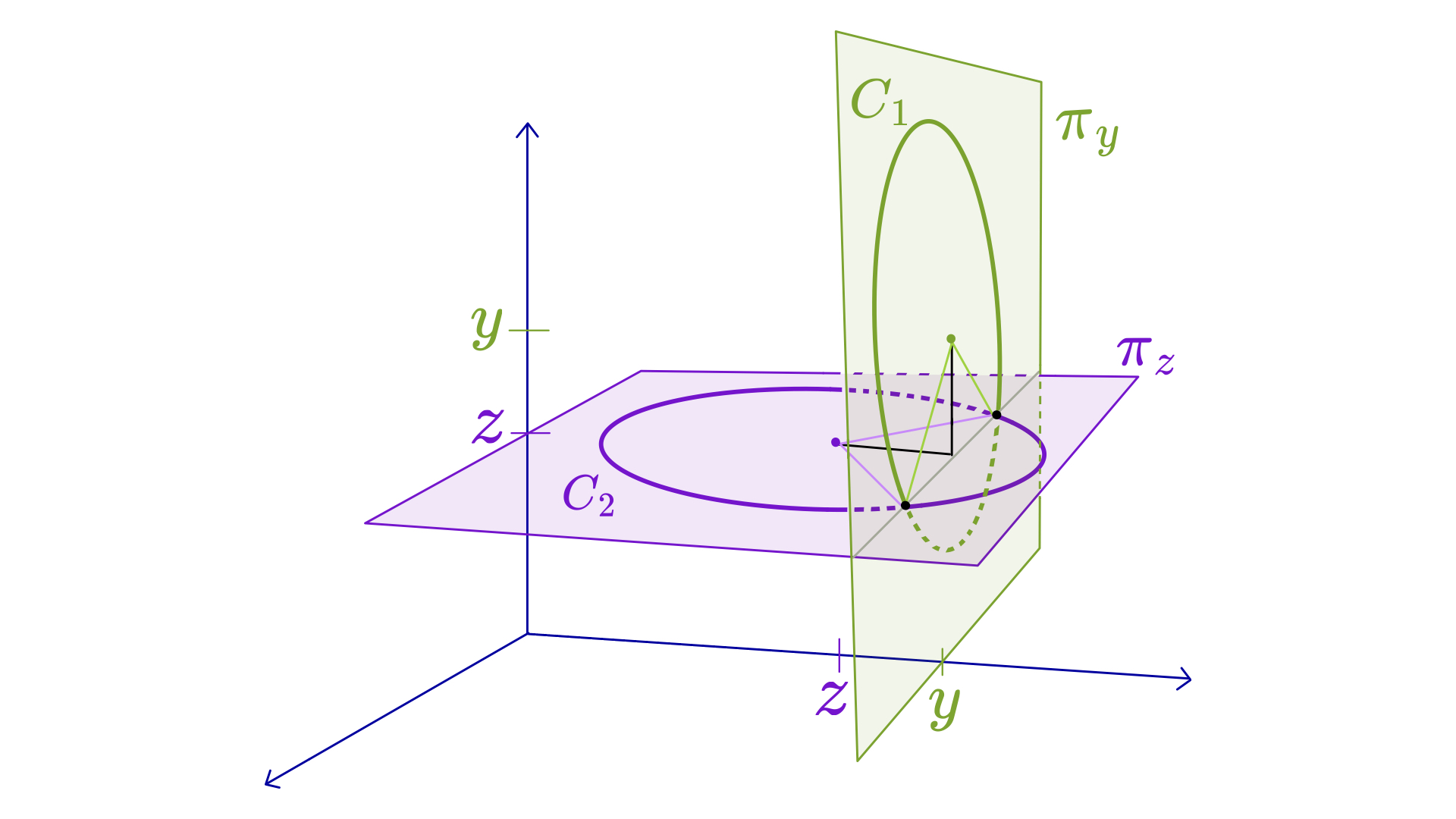}
	\caption{$C_1$ and $C_2$ intersect in the points $\left(y,z,\pm \sqrt{1-(z-y)^2}\right)$.}
	\label{fig: the-order-is-needed-puc}
\end{figure}

It would be tempting to try to consider $\mathbf{P_C}$ with the order given just by reverse inclusion, since this is enough for the cases of Mazurkiewicz sets and Hamel bases (see Section \ref{section: appendix}).
However, using this forcing, amalgamation does not not work. 
Consider $x,y,z\in \RR$ and $p,p_1,p_2$ as in Equations \ref{eq: amalg puc p in V[x]}--\ref{eq: amalg puc p2 in V[x,z]} in the proof of Lemma \ref{lemma: PUC amalgamation}. 
Let $\pi_y$ and $\pi_z$ be the planes that consist of all the points in $\RR^3$ with first coordinate $y$ and second coordinate $z$ respectively.
Assume $|z-y|<1$, and that $p_1$ contains a circle $C_1$ and $p_2$ contains a circle $C_2$ described as follows: 
\begin{align*}
	C_1 &\text{ is the only unit circle contained in } \pi_y \text{ and origin } (y,y,0), \text{ and}\\
	C_2 &\text{ is the only unit circle contained in } \pi_z \text{ and origin } (z,z,0). 
\end{align*}
Figure \ref{fig: the-order-is-needed-puc} shows that the circles $C_1$ and $C_2$ will intersect in the points 
\[\left(y,z, \sqrt{1-(z-y)^2}\right), \left(y,z, -\sqrt{1-(z-y)^2}\right).\]
This contradicts amalgamation for $(\mathbf{P_C}, \supseteq)$.
\\

Now we are ready to prove the main theorem of this section. 

\begin{Teo}[Corollary of Theorem \ref{th: partition no wo}] \label{th: PUC - no wo}
	Let $\mathbf{Q}$ be the finite support product of $\omega_1$-many copies of  Cohen forcing, let $g$ be a $\mathbf{Q}$-generic filter over $V$. 
	Let $\mathbf{P_C}$ be the forcing poset in $V[g]$ described in Definition \ref{def: forcing PUC}. 
	Let $h$ be a $\mathbf{P_C}$-generic filter over $V[g]$, and let $\mathcal{P}=\cup h$.
	Then 
	\[L(\RR, \mathcal{P})^{V[g,h]}  \models \mathsf{ZF} + \mathsf{DC} +  \neg \mathsf{WO}(\RR) + \mathcal{P}\text{ is a partition of } \RR^3 \text{ into unit circles}. \]
\end{Teo} 

\begin{Dem}
	We will apply Theorem \ref{th: partition no wo}. 
We have shown that $\mathbf{P_C}$ adds a real partition.
	It is not a trivial forcing because for any $x\in \RR^{V[g]}$, $V[x]$ is a model of $\mathsf{AC}$, and thus has a PUC (see Theorem \ref{th: zfc PUC}).
	Lemmas \ref{lemma: PUC extendability}, \ref{lemma: PUC amalgamation} and \ref{lemma: PUC forcing sigma closed} show, respectively, that $\mathbf{P}$ satisfies extendability, amalgamation, and is $\sigma$-closed.
	We can then apply Theorem \ref{th: partition no wo} and obtain the desired conclusion. 
\end{Dem}

\section{A PUC in the Cohen-Halpern-Lévy model}	\label{section: PUC in the Cohen model}

We follow the structure of the proof in \cite{Beriashvili2018}, which shows that the Cohen model $H$ has is a Hamel basis of $\mathbb{R}$.
For this, we will need a stronger version of the notion of amalgamation (see Definition \ref{def: omega amalagamation}) proven to be valid for PUCs in Lemma \ref{lemma: PUC amalgamation} and therefore a stronger version of the Lemma \ref{lemma: trascendence degree one cohen} to be able to prove it, which is
\cite[Theorem 3.13]{Fatalini2024}. \\

\begin{Teo}\label{th: PUC in Cohen model}
	Let $\mathbf{C}(\omega)$ denote the finite support of $\omega$-many copies of $\mathbf{C}$, let $g$ be a $\mathbf{C}(\omega)$-generic filter over $L$ and $A$ be the set of Cohen reals added by $g$. 
	Let $H$ be the Cohen-Halpern-Lévy model as described in  Definition \ref{def: prereq cohen model}.
	Then 
	\[H = \mathsf{HOD}_A^{L[g]} \models \text{There is a PUC } + \neg \mathsf{AC}_\omega. \]
\end{Teo}

\begin{Dem}
	Using Theorem \ref{th: in H A has no countable subset}, we deduce that $\mathsf{AC}_\omega$ does not hold in $H$.
	So we only need to prove that there is a partition of unit circles in $H$. 
	
	Work inside $H$. 
	We will construct a family $\{p_Y\}_{Y\in [A]^{<\omega}}$
	so that each $p_Y$ is a partition of unit circles in $L[Y]$, and
	for each $Y\in [A]^{<\omega}$ such that $Y\subseteq X \in [A]^{<\omega}$, then $p_X \leq p_Y$, where $\leq$ is defined as $\leq_{\mathbf{P_C}}$ in Definition \ref{def: forcing PUC}.
	We will do so recursively on $n=|Y|$.
	
	For $n=0$: notice that $L$ has a PUC by Theorem \ref{th: zfc PUC}. 
	Let $p_\emptyset\in L$ be the $<_\emptyset$-least PUC in $L$.
	We do so using the global well order $<_\emptyset$ (see Lemma \ref{lemma: prereq H has the well orders of L[a]}).
	Suppose we already defined $p_Y$ for all $Y\subseteq A$ with $|Y|\leq n$. 
	
	Let $X$ be a subset of $A$ of size $n+1$. 
	Consider 
	\[p^\ast = \bigcup_{Y\subsetneq X} p_Y.\]
	Let $Y,Y'\subsetneq X$ and $Y\neq Y'$. We claim that $p_Y$ and $p_{Y'}$ are compatible, namely, its union is a family of disjoint unit circles. 
	If either $Y$ or $Y'$ is a subset of the other, for example, $Y\subseteq Y'$, then by inductive hypothesis $p_{Y'}\leq p_Y$, and hence $p_{Y'}\cup p_Y=p_{Y'}$.
	If not, then consider $Z=Y\cap Y'$. 
	Then $Z\subseteq Y, Y'$. 
	Recall that $\mathbf{C}\cong \mathbf{C}^k$ for any $k<\omega$ (Theorem \ref{th: prereq cohen is the only countable forcing}).
	By the proof of Lemma \ref{lemma: PUC amalgamation}, we get that $p_Y$ and $p_Y'$ are compatible. 
	To check if $p^\ast$ is a family of disjoint unit circles, we need to take two circles, and check whether they intersect. 
	By the pairwise compatibility of $\{p_Y \mid Y \subsetneq X \}$ and recalling that intersection between two circles is absolute, we get that $p^{\ast}$ is a family of disjoint unit circles. 
	
	We need to prove that in $L[X]$ there is a partition of unit circles $p_X$ such that $p_X\leq p_Y$ for all $Y\subsetneq X$. 
	In particular, we need $p^\ast\subseteq p_X$. 
	We will proceed in a way similar to the proof of Lemma \ref{lemma: PUC extendability}.
	
	Work in $L[X]$.
	Let $\{x_\alpha\}_{\alpha<\mathfrak{c}}$ be an enumeration of the points in $\RR^3\backslash \cup p^\ast$.
	Here $\cup p^\ast$ is the union of all the circles given by $p^\ast$ as computed in $L[X]$.
	We will recursively define $p_\alpha$ for $\alpha<\mathfrak{c}$. 
	For $\alpha=0$, set $p_0=p^\ast$. 
	Suppose that $p_\beta$ is defined for all $\beta< \alpha$. 
	If $\alpha$ is a successor ordinal of the form $\beta +1$, and $x_\beta \in \cup p_\beta$ (namely, $x_\beta$ is covered by a circle in $p_\beta$), take $p_{\beta+1}=p_\beta$. If $x_\beta \not \in \cup p_\beta$, we will pick a unit circle $C_{\beta}$ such that $x_\beta \in C_{\beta}$ and $C_{\beta} \cap C= \emptyset$ for all $C \in p_{\beta}$. 
	Assuming we can choose such a $C_\beta$, we define $p_{\beta +1}=p_{\beta}\cup \{C_{\beta}\}$.
	Finally, if $\alpha$ is a limit ordinal, define $p_\alpha= \bigcup_{\beta <\alpha} p_{\beta}$. 
	We need to check that the construction is possible, namely, that we can choose such a circle $C_\beta$.
	
	Since $C_\beta$ will have radius $1$, we only need to choose a center $o_{\beta}$ of the circle and a vector $n_\beta$ normal to the plane in which $C_\beta$ will be contained in. 
	We want to choose $n_\beta$ to be different from all the normal vectors of circles in $p_\beta$.
	Let
	\[\RR^\ast= \bigcup_{Y\subsetneq X} (\RR\cap L[Y]).\]
	Notice that if $C\in p^{\ast}$, then its parameters $(o,n)$ would be in $\RR^\ast$. 
	Since $|\RR \backslash \RR^\ast|=\mathfrak{c}$,
	we have continuum many options for $n_\beta$ that are different from all the normal vectors of circles in $p^{\ast}$. 
	Since $p_\beta=p^\ast\dot{\cup} \tilde{p}_\beta$ and $|\tilde{p}_\beta|\leq |\beta|$, we have to also avoid choosing at most $|\beta|$-many normal vectors (one per circle in $\tilde{p}_\beta$).
	Since $|\beta|<\mathfrak{c}$, we can choose $n_\beta$ with the desired property. 
	This implies that the plane $\pi_\beta$ (determined by $n_\beta$ and $x_\beta$) in which $C_\beta$ will be contained is different from all the planes containing circles in $p_\beta$.
	Therefore, $|\pi_\beta \cap C|\leq 2$ for every circle $C\in p_\beta$.
	
	Additionally, it is clear that we have to choose $o_\beta$ in $\pi_{\beta}$ so that the distance between $x_\beta$ and $o_\beta$ is equal to $1 $.
	The locus of such a point is a circle $C$ contained in $\pi_\beta$ with center $x_\beta$ and of radius $1$.
	
	Let us choose $o_\beta\in \pi_\beta$ such that $o_\beta\not \in \overline{F_\beta}$ (i.e. at least one coordinate is not an element of $\overline{F_\beta}$), where 
	\[ F_\beta= \text{ the minimal field containing }\mathbb{R}^\ast \cup
	\coor(\tilde{p}_\beta) \cup
	\coor(\pi_{\beta}, x_\beta).\]
	
	Recall that $\coor(\tilde{p}_\beta)=\bigcup_{\delta<\beta}\coor(o_\delta, \pi_{\delta})$, and therefore it has cardinality at most $|\beta|<\mathfrak{c}$. 
	Applying \cite[Theorem 3.13]{Fatalini2024} to this context, we know that $\RR=\RR^{L[X]}$ has transcendence degree $\mathfrak{c}$ over the minimal field containing $\RR^\ast$. 
	Also, $\coor(\tilde{p}_\beta) \cup
	\coor(\pi_{\beta}, x_\beta)$ has cardinality $|\beta|<\mathfrak{c}$.
	So we can conclude $|\RR \backslash \overline{F_\beta}|= \mathfrak{c}$. 
	Finally, we can choose $o_\beta$ such that $o_\beta \not\in \overline{F_\beta}$ because
	we still have one degree of freedom after prescribing \(o_{\beta} \in \pi_{\beta}\) and \(d (o_{\beta}, x_{\beta})=1\). 
	
	Let $C_\beta$ be the circle determined by $(o_\beta, n_\beta)$. 
	We have to check that it satisfies the required properties for the recursive construction. 
	Clearly $x_\beta \in C_\beta$. 
	Now fix $\tilde{C}\in p_\beta$. 
	We want to show $C_\beta \cap \tilde{C}=\emptyset$. 
	Suppose there is some \(t\in C_{\beta}\cap \tilde{C}\).
	If $\tilde{C}\in p^\ast$, its parameters $(o, \pi)$ belong to \(\mathbb{R}^\ast \).
	We can calculate $o_{\beta}$ from $(o, \pi, x_{\beta}, \pi_{\beta})$ ``algebraically'': $\tilde{C}$ can be computed from $(o, \pi)$, $t$ can be computed from $(\tilde{C}, \pi_\beta)$, and $o_\beta$ can be computed from $(t, x_\beta, \pi_\beta)$ (see Figure \ref{fig: PUC-fig3}).
	This means all the coordinates of $o_\beta$ belong to $\overline{F_\beta}$, contradicting the choice of $o_{\beta}$. 
	
	The case in which $\tilde{C}\in \tilde{p}_\beta$ is analogous. 
	$\tilde{C}$ must have been added in some step $\delta +1<\beta$.
	We can then calculate $o_{\beta}$ from $(o_\delta, \pi_\delta, x_{\beta}, \pi_{\beta})$ ``algebraically'' in the same fashion, from which we get the same contradiction. 
	
	Take $\tilde{p}=\bigcup_{\alpha<\mathfrak{c}} p_\alpha$.
	For any two circles $C, D \in \tilde{p}$, we want to show that $C\cap D=\emptyset$.
	This is clear for $C,D$ added in step $0$, namely, $C,D\in p^\ast$. 
	If they were added in different steps, for example, $D$ strictly after $C$, there is $\alpha<\omega_1$ such that $D=C_\alpha$ and $C\in p_\alpha$.
	By construction, $C_\alpha\cap C =\emptyset$, so we can conclude $\tilde{p}$ is a family of disjoint unit circles. 
	Moreover, for any $r\in \RR^3$, either $r\in \cup p^\ast$ or there is	
	$\alpha < \mathfrak{c}$ such that $r= x_\alpha$. 
	In the first case, $r\in \cup \tilde{p}$, since $p^\ast=p_0\subseteq \tilde{p}$.
	In the second case, $r\in \cup p_{\alpha+1} \subseteq \cup \tilde{p}$ by construction. 
	Thus, $\tilde{p}$ is a partition of $\RR^3$ into unit circles that extends $p^\ast$.
	
	Moreover, $\tilde{p}\leq p_Y$ for every $Y\supsetneq X$: Clearly, $\tilde{p}\supseteq p_Y$. Also, 
	$Y\in L[X]$.
	Fix $C\in \tilde{p} \backslash p_Y$. 
	By construction, $o\not \in \overline{{\RR^\ast}(\coor(\pi))}$ and $\pi \not\in \RR^\ast$.
	Since $\RR^\ast \supseteq \RR^{L[Y]}$, then $o\not \in \overline{\RR^{L[Y]} (\coor(\pi))}$ and $\pi \not \in \RR^{L[y]}$.
	Therefore $\tilde{p}\leq p_Y$ for every $Y\supsetneq X$ as we wanted. \\
		
		We have just proved that in $L[X]$ there is a partition of unit circles that extends $p^\ast$ and that is below $p_Y$ (according to $\leq_{\mathbf{P_C}}$) for each $Y\subsetneq X$. 
		In $H$, let $p_X$ be the $<_{X}$-least such partition (see Lemma \ref{lemma: prereq H has the well orders of L[a]}).
		Finally, define $p=\bigcup_{Y\in[A]^{<\omega}} p_Y$. 
		We claim $p$ is a partition of unit circles (in $H$): 
		Clearly, it is a family of unit circles. If $C_0\in p_X$ and $C_1\in p_Y$, then $C_0, C_1 \in p_{X\cup Y}$ which is a partition of unit circles in $L[X\cup Y]$; therefore, $C_0$ and $C_1$ are disjoint.
		Let $r\in \RR^3$. 
		By Theorem \ref{th: in H A has no countable subset} there is some $Y\in [A]^{<\omega}$ such that $r\in \RR^3 \cap L[Y]$. Therefore $r$ is covered by some circle $C\in p_Y$. 
	\end{Dem}
	
	We will capture the main obstacle in the proof of Theorem \ref{th: PUC in Cohen model} by the following definition. 
	
	\begin{Def} \label{def: omega amalagamation}
		Let $g$ be a $\mathbf{C}(\omega)$-generic filter over $L$, and let $A$ be the set of reals added by $g$. 
		Assume $\mathbf{P}\in L[g]$ is a forcing that adds a real partition according to $\psi$ as in Definition \ref{def: forcing adds a real partition}. 
		
		Let $X$ be a finite subset of $A$. 
		We say that $\{p_Y\}_{Y\subsetneq X}\subseteq \mathbf{P}$ is a \concept{compatible family of conditions} iff for every $Z\subseteq Y\subsetneq X$ we have that 
		\[L[Y]\models \psi(p_Y) \text{ and } p_Y\leq_{P} p_Z.\]
		
		We say that $\mathbf{P}$ satisfies \concept{($<\omega$)-amalgamation} if 
		for all $X\in [A]^{<\omega}$ and for all family of compatible conditions $\{p_Y\}_{Y\subsetneq X}$ in $\mathbf{P}$,
		we have that 
		\[L[X]\models \bigcup_{Y\subsetneq X} p_Y \text{ is a partial condition (as in Def. \ref{def: partial condition and extendability})}\]
		and moreover, there is $p_X\in L[X]$ such that 
		\[L[X] \models \psi(p_X) \text{ and } p_X \leq_{\mathbf{P}} p_Y \text{ for every } Y\subseteq X.\]   
	\end{Def}
	
	\begin{Remark}
		The proof of Theorem \ref{th: PUC in Cohen model} actually showed that $\mathbf{P_C}$ satisfies $(<\omega)$-amalgamation in $L[g]$, where $g$ is $\mathbf{C}(\omega)$-generic over $L$.
		The proof of the existence of a Hamel basis in $H$ \cite[Theorem 2.1]{Beriashvili2018} essentially shows that the partial order $\mathbf{P_H}$ defined in Corollary \ref{coro: hamel basis - no wo} satisfies  ($<\omega$)-amalgamation in $L[g]$ as well.
		This is a strategy that has been proven to work to find some paradoxical sets in $H$. 
		However, it is hard to see whether the forcing $\mathbf{P_M}$ for Mazurkiewicz sets satisfies ($<\omega$)-amalgamation. 
		It has been shown that the Cohen model $H$ contains a Mazurkiewicz set \cite[Corollary 0.3]{Beriashvili2022}, but using other construction.
	\end{Remark}

\newpage
\section{Appendix} \label{section: appendix}

To show the usefulness of the theorems in Section \ref{section: gral setup} (specially Theorem \ref{th: general set up}), we will apply them in this section to two other paradoxical sets: Hamel bases and Mazurkiewicz sets, recovering known results in the literature. 

\subsection{Application 2: Hamel bases}\label{subsection: hamel bases}

In this subsection we will show that we can apply the methods of Section \ref{section: gral setup} to this paradoxical set, getting a models of $\mathsf{ZF}+\mathsf{DC}+\neg \mathsf{WO}(\RR)$ with a Hamel basis.
We can apply Theorem \ref{th: partition no wo} to the case of Hamel bases recovering the result   \cite[Theorem 1.1]{Schindler2018}. 
Moreover, to change the Cohen reals for Sacks reals in the forcing $\mathbf{Q}$ and get back \cite[Theorem 5.1]{Brendle2018}, we will apply Theorem \ref{th: general set up}.

\begin{Def}
	Let $H\subseteq \RR$. 
	We say that $H$ is a Hamel basis if it is a basis of $\RR$ as a vector space over $\mathbb{Q}$, namely, 
	a maximal linearly independent set over $\mathbb{Q}$.
\end{Def}

It is clear that $\mathsf{ZF}+\mathsf{WO}(\RR)$ implies there is a Hamel basis: one can construct a Hamel basis by extending recursively a linearly independent set until it is maximal, always adding the first real which is not in the span of the linearly independent set taken so far. 

\begin{Coro}[of Theorem \ref{th: partition no wo}] \label{coro: hamel basis - no wo}
	Let $\mathbf{Q}$ be the finite support product of $\omega_1$-many copies of Cohen forcing, and let $g$ be a $\mathbf{Q}$-generic filter over $V$. 
	Let $\mathbf{P_H}$ be the forcing poset in $V[g]$ given by
	\[p \in \mathbf{P_H} \iff \exists x \in \RR:\; V[x]\models p \text{ is a Hamel basis}, \]
	ordered by reverse inclusion. 
	Let $h$ be a $\mathbf{P_H}$-generic filter over $V[g]$, and let $\mathcal{P}=\cup h$.
	Then 
	\[L(\RR, \mathcal{P})^{V[g,h]}  \models \mathsf{ZF} + \mathsf{DC} + \neg \mathsf{WO}(\RR)+\mathcal{P} \text{ is a Hamel basis}. \]
\end{Coro} 

\begin{Dem}
	We want to apply Theorem \ref{th: partition no wo}, so let us verify its hypotheses. 
	
	\begin{enumerate}
		\item \emph{$\mathbf{P_H}$ adds a real partition:} The natural way to define a Hamel basis works for the role of $\psi$. Then, if $(x,p)$ is such that 
		\[V[x]\models p \text{ is a Hamel basis},\]
		with $x\in \RR^{V[g]}$, then by definition there is in $V[x]$ a finite sequence $\vec{s}$ of elements of $p$ and a finite sequence  $\vec{q}$ of rational numbers such that $\vec{s}\cdot\vec{q}=x$.
		\item \emph{$\mathbf{P_H}$ is $\sigma$-closed:} This follows from Lemma \ref{lemma: sigma close for partitions}.
		\item \emph{$\mathbf{P_H}$ satisfies extendability:} 
		The partial conditions $p$ of $\mathbf{P_H}$ are $\mathbb{Q}$-linearly independent subsets of reals in some $V[x]$ for $x\in \RR^{V[g]}$. 
		Because a Hamel basis is a maximal linearly independent set, then already in $V[x]$ there is a Hamel basis $\bar{p}$ such that $\bar{p}\supseteq p$. 
		Notice that $p\in \mathbf{P_H}$ and, in case $p$ was a condition, $\bar{p}\leq_{\mathbf{P_H}} p$ trivially.
		
		\item  \emph{$\mathbf{P_H}$ satisfies amalgamation:} Work in $V[g]$ and fix $p\in \mathbf{P_H}$, and $g_1$ and $g_2$ mutually $\mathbf{Q}$-generic over $V[p]$. 
		Fix $p_1\in \mathbf{P_H}\cap V[p,g_1]$ and $p_2\in \mathbf{P_H}\cap V[p,g_2]$ such that $p_1\leq_{\mathbf{P_H}}$ and $p_2\leq_{\mathbf{P_H}} p$. 
		We want to show $p_1$ and $p_2$ are compatible conditions, i.e., we want to show that
		$\bar{p}=p_1\cup p_2$ is linearly independent. 
		
		This is not a new argument (see \cite[Claim 3]{Brendle2018}), but we reproduce it here for completeness. 
		Suppose $\bar{p}$ is not linearly independent over $\mathbb{Q}$. 
		Then there is $\vec{s}\in [\bar{p}]^{<\omega}$, $ \vec{q}\in [\QQ]^{<\omega}$ such that $\vec{s} \cdot \vec{q}=0$. 
		Separating terms accordingly, we can write 
		\[   \vec{s_0} \cdot \vec{q_0} +\vec{s_1} \cdot \vec{q_1}+ \vec{s_2} \cdot \vec{q_2}=0, \]
		where $\vec{s_0}\subseteq p$, $\vec{s_1} \subseteq p_1\backslash p, \vec{s_2} \subseteq p_2 \backslash p$; $s=s_0\cup s_1 \cup s_2$ and $q=q_0\cup q_1 \cup q_2$.
		Now, notice that
		\[\vec{s_1} \cdot \vec{q_1}=-\vec{s_0} \cdot \vec{q_0} - \vec{s_2} \cdot \vec{q_2} \in \RR^{V[p,g_1]}\cap \RR^{V[p,g_2]}=\RR^{V[p]}. \]
		Then, $ - \vec{s_1} \cdot \vec{q_1}$ is a real number in $V[p]$.
		Since $p$ is a Hamel basis there, there are $\vec{t_0}\in [p]^{<\omega}$ and $\vec{r_0}\in [\QQ]^{<\omega}$ such that 
		\[ \vec{t_0} \cdot \vec{r_0} =  - \vec{s_1} \cdot \vec{q_1}. \]
		Then we get
		\[ \vec{t_0} \cdot \vec{r_0} + \vec{s_1} \cdot \vec{q_1} =0. \]
		Notice that $s_1\subseteq \RR^{V[p,g_1]}\backslash \RR^{V[p]}$, so $t_0\cap s_1 = \emptyset$.
		Since $p_1$ is linearly independent, we get that $\vec{r_0}=\vec{0}$ and $\vec{q_1}=\vec{0}$. 
		Coming back to the first equation, we have 
		\[\vec{s_0} \cdot \vec{q_0} + \vec{s_2} \cdot \vec{q_2}=0.\]
		Similarly, we obtain $\vec{q_0}=\vec{0}$, and $\vec{q_2}=\vec{0}$. 
		Therefore $\vec{q}=\vec{0}$, as we wanted. \\
		Now, let $r_1$ and $r_2$ be reals such that $p_1$ and $p_2$ are Hamel bases in $V[r_1]$ and $V[r_2]$ respectively.
		In $V[r_1\oplus r_2]$, we can extend $p_1\cup p_2$ to a Hamel basis $p^{\ast}$. 
		Then $p^{\ast}$ witnesses the compatibility of $p_1$ and $p_2$ in $V[g]$. 
	\end{enumerate}
	Applying Theorem \ref{th: partition no wo}, we get that 
		\[L(\RR, \mathcal{P})^{V[g,h]}  \models \mathsf{ZF} + \mathsf{DC} + \neg \mathsf{WO}(\RR)+\mathcal{P} \text{ is a Hamel basis}. \]
\end{Dem}

\begin{Coro}[of Theorem \ref{th: general set up}] \label{coro: hamel basis - sacks forcing}
	Let $\mathbf{Q}$ be the countable support product of $\omega_1$-many copies of Sacks forcing, and let $g$ be a $\mathbf{Q}$-generic filter over $V$. 
	Let $\mathbf{P}$ be the forcing poset in $V[g]$ given by
	\[p \in \mathbf{P} \iff \exists x \in \RR: \; V[x]\models p \text{ is a Hamel basis of } \RR, \]
	ordered by reverse inclusion. 
	Let $h$ be a $\mathbf{P}$-generic filter over $V[g]$, and let $\mathcal{P}=\cup h$. 
	Then 
	\[L(\RR, \mathcal{P})^{V[g,h]}  \models \mathsf{ZF} + \mathsf{DC} + \neg\mathsf{WO}(\RR)+\mathcal{P} \text{ is a Hamel basis}. \]
\end{Coro} 

	\begin{Dem}
$\mathbf{Q}$ is homogeneous and $\mathbf{Q}\times \mathbf{Q} \cong \mathbf{Q}$, see \cite[Claim 2]{Brendle2018}.
$\mathbf{P}$ is real absolute and $\mathbf{P}$ and $\mathbf{Q}$ are real-alternating by the same argument in the proof of Theorem \ref{th: partition no wo}. Also, $\mathbf{P}$ is $\sigma$-closed by the same argument of the proof of Lemma \ref{lemma: sigma close for partitions}. 
We only need to prove that $\mathbf{P}$ is  a $\mathbf{Q}$-balanced forcing over $V$. Notice that proving amalgamation is enough and that in the proof of amalgamation in Corollary
\ref{coro: hamel basis - no wo} we did not use any property of Cohen reals, so we obtain amalgamation in this case as well. 
\end{Dem}

\subsection{Application 3:  Mazurkiewicz sets}\label{subsection: maz sets}

In this section we will deal with another example of paradoxical set with geometrical flavor. Mazurkiewicz proved in 1914 that these particular sets existed \cite{Mazurkiewicz1914}, using the Axiom of Choice. 
Actually, a well-ordering of the reals is enough to carry on that proof.
We will apply Theorem \ref{th: partition no wo} to this case, and therefore recover a result of existence of a model of $\mathsf{ZF}+\mathsf{DC}+\neg \mathsf{WO}(\RR)$ with a Mazurkiewicz set (done in unpublished notes of Beriashvili and Schindler \cite{Beriashvili2019} and also mentioned in \cite{Brendle2018}). 

\begin{Def}
	Let $M\subseteq \RR^2$. 
	We say that $M$ is a \concept{Mazurkiewicz set} (also called \concept{two point set}) if for every line $l$ in $\RR ^2$, $| M\cap l | = 2$. 
\end{Def}

Unlike the example of Hamel bases and similarly to the case of PUCs, partial two-point sets may not be extendable to complete two-point sets.
This makes that the conditions of \emph{extendability} and \emph{amalgamation} that Theorem \ref{th: partition no wo} requires are harder to get, so they will be treated in separate lemmas, Lemma \ref{lemma: mazurkiewicz extendability} and Lemma \ref{lemma: mazurkiewicz amalgamation} respectively.

\begin{Def}\label{def: forcing Mazurkiewicz}
	Let $V$ be a model of $\mathsf{ZFC}$. 
	Let $\mathbf{Q}$ be the finite support product of $\omega_1$-many copies of Cohen forcing, let $g$ be a $\mathbf{Q}$-generic filter over $V$. 
	Let us define $\mathbf{P_{M}}$ as the forcing poset in $V[g]$ given by
	\[p \in \mathbf{P_M} \iff \exists x \in \RR\; V[x]\models p \text{ is a Mazurkiewicz set}, \]
	ordered by reverse inclusion. 
\end{Def}

Observe that $\mathbf{P_M}$ adds a real partition. 
For this, notice that $p$ is a two-point set iff:
\[p\subseteq \RR^2  \wedge \left(\forall s \in [p]^3 \psi_1(s) \right) \wedge \left( \forall r\in\RR^3 \exists s \in [p]^{2} \psi_2(r,s)\right), \]
where $\psi_1(s)$ iff ``the elements of $s$ are not collinear'' and $\psi_2(r,s)$  iff ``the elements of $s$ belong to the line given by $r$''.
Clearly, $\psi_1$ and $\psi_2$ are $\bm{\Delta_0}$ and therefore absolute.
Also, for any pair $(x,p)$ as before, since $p$ is a Mazurkiewicz set in $V[x]$ and $x$ is a real, there is $s\in [p]^2$ such that $s$ is contained in the line $l_x = \{ (x,y) \mid y \in \RR \}$. 
Then $x$ can be computed from $s$, by taking the first coordinate of any of its elements. 
\\

\begin{Remark}
	To show that $\mathbf{P_M}$ adds a real partition we implicitly assumed that we have fixed a representation of the lines in $\RR^2$ by points in $\RR^3$. 
	For example, let
	\[S= \{(a,b,c)\in \RR^3 \mid c=1 \lor (c=0 \land a=1)  \}. \]
	Then for any $(a,b,c)\in S$ we can define a line $l$ in $\RR^2$ by 
	\[l= \{(x,y)\in \RR^2 \mid ax+b=cy \}. \]
	Conversely, for any line $l$ there is a \emph{unique} set of parameters $(a,b,c)\in S$ that determine $l$ in this way.
	Formally, we define $\psi_2(r,s)$ so that also holds true in any case that $r$ does not belong to the image of such representation.
\end{Remark}\\

\begin{Not}
	Let $p\subseteq \RR^2$. 
	By $\tup{p}$ we denote the set
	\[\{l \text{ line } \mid \exists s_1, s_2 \in p \text{ and } l = l(s_1, s_2) \}. \]
	We will frequently consider the set $\cup \tup{p}$. 
	Notice that $\tup{p}$ is a set of lines and $\cup \tup{p}$ is instead a set of points in $\RR^2$.
\end{Not}

\begin{Not}
	For any line $l$, we will confuse it (the geometrical object) with its representation as an element of the set $S$ described above. 
	Similarly to the case of PUCs, we will consider the \concept{parameters} or \concept{coordinates} of $l$ as the set $\coor(l)=\{a,b,c\}$, where $a,b,c$ are such that $(a,b,c)\in S$ and 
	$l= \{(x,y)\in \RR^2 \mid ax+b=cy \}$. 
	Moreover, for any model $M$ we will write ``$l\in M$'' as a short form of ``$\coor(l)\in M$''.
	
	Similarly, if $r\in \RR^n$ with $r=(r_0,\dots,r_{n-1})$, we write $\coor(r)$ to denote the set $\{r_0,\dots, r_{n-1}\}$. 
	Furthermore, if $R\subseteq \RR^n$ or $R$ is a set of lines, we denote the set $\bigcup \{\coor(r) \mid r\in R \}$ by $\coor(R)$.
\end{Not}

\begin{Lema}\label{lemma: mazurkiewicz extendability}
	Let $\mathbf{Q}$ be the finite support product of $\omega_1$-many copies of Cohen forcing, let $g$ be a $\mathbf{Q}$-generic filter over $V$. 
	Then $\mathbf{P_{M}}$ in $V[g]$ satisfies extendability. 
\end{Lema}

\begin{Dem}
	Looking at Definition \ref{def: partial condition and extendability}, we need to prove that for any partial condition $p$ in $V[x]$ ($x\in\RR^{V[g]} $), there is a condition $\bar{p}$ and a real $\bar{x}$ such that $\bar{p}\supseteq p$, and $x \in V[\bar{x}]$. 
	Recall that $\leq_{\mathbf{P_M}}=\supseteq \restriction \mathbf{P_M}$.
	So, if $p$ is a condition, then $\bar{p} \leq_{\mathbf{P_M}} p$.
	
	Fix $x\in \RR^V[g]$, and let $p$ be a partial condition in $V[x]$. 
		Notice that $p$ is a subset of $\RR^2 \cap V[x]$ such that no three points are collinear.
	Let $\gamma<\omega_1$ be such that $x\in V[g\restriction\gamma]$, let $y$ be $\bigcup (g\restriction\{\gamma\})$, and define $\bar{x}$ as $x\oplus y$.
	We will use a variation of the proof of existence of Mazurkiewicz sets in $\mathsf{ZFC}$ in order to construct a Mazurkiewicz set $\bar{p}$ inside $V[\bar{x}]$ extending $p$. 
	
	Work inside $V[\bar{x}]$. 
	By absoluteness, no three points in $p$ are collinear. 
	Notice that $\tup{p} \subseteq  \{l \text{ line } \mid l\in V[x] \}$.
	
	Let $\{l_\alpha\}_{\alpha<\mathfrak{c}}$ be an enumeration of  all the lines excepting the ones in $\tup{p}$.
	We will recursively define $p_\alpha \subseteq \RR^2$ for $\alpha<\mathfrak{c}$. 
	For $\alpha=0$, take $p_0=p$.
	Now suppose $p_{\beta}$ is defined for all $\beta < \alpha$. 
	If $\alpha$ is a successor ordinal, namely $\alpha=\beta +1$, we will take $r\subseteq l_\beta$ such that $|(p_\beta \cup r) \cap l_{\beta}|=2$, and such that for each element of $r$, there is a coordinate of it that is not in the real algebraic closure of $F_\beta$, where 
	\[ F_\beta = \text{the minimal field containing } (\RR\cap V[x]) \cup \coor(p_\beta) \cup  \coor(l_\beta).\]

	This implies that no element of $r$ is in $l_\beta \cup \tup{p_\beta}$ since any intersection point $l\cap l_\beta$ with $l\in \tup{p_\beta}$ would have both coordinates in $\overline{F_\beta}$.
	Define $p_\alpha=p_\beta\cup r $.
	If $\alpha$ is a limit ordinal, take $p_\alpha=\bigcup_{\beta<\alpha} p_\beta$.
	
	We have to check that the construction is possible, namely that such $r$ exists. 
	First, we will show that $|p_\beta \cap l_{\beta}|\leq 2$ for all $\beta<\mathfrak{c}$. 
	Suppose $\beta$ is the first ordinal such that $|p_\beta \cap l_{\beta}|\geq 3$. 
	Let $x, y, z$ be three points in $p_\beta \cap l_{\beta}$, named alphabetically by the order of being added to the construction. 
	If $x\in p$ we say $x$ was added in the step 0. 
	Since $p$ is a partial condition, $z\not \in p$ and $z$ should have been added at some step $\delta+1$ which is of course different from 0, which means $z\in l_\delta$.
	By construction, $z\not \in \cup \tup{p_\delta}$. 
	This is a contradiction, since $l(x,y)\in \tup{p_\delta}$. 
	Thus, $|p_\beta \cap l_{\beta}|\leq 2$ for all $\beta<\mathfrak{c}$.
	
	The rest of the proof consists of showing that $\RR \backslash \overline{F_\beta}$ has at least two points so that we can choose $r$.  
	Notice that $p_\beta=p \dot{\cup} \tilde{p}_\beta$, where $|\tilde{p}_\beta|\leq |\beta|<\mathfrak{c}$. 
	Since $p\subseteq V[x]$, we can write 
	\[F_\beta = \text{the minimal field containing } (\RR\cap V[x]) \cup \coor(\tilde{p}_\beta) \cup  \coor(l_\beta). \]
	
	Thus, $F_\beta= \RR^{V[x]}(S)$, 
where $S$ is a set of cardinality strictly less than $\mathfrak{c}$.
	Applying Lemma \ref{lemma: trascendence degree one cohen} and recalling that $y$ was a Cohen real over $V[x]$, we know that the transcendence degree of $\RR=\RR^{V[\bar{x}]}$ over $\RR\cap V[x]$ is $\mathfrak{c}$.
	Therefore $\RR \backslash \overline{F_\beta}$ is actually of cardinality $\mathfrak{c}$, and there are enough possibilities to choose $r$ from.
	
	Take $\bar{p}=\bigcup_{\alpha<\mathfrak{c}} p_{\alpha}$. 
	Then $\bar{p}$ is a Mazurkiewicz set in $V[\bar{x}]$ and it contains $p$.
\end{Dem}

In the proofs of Lemma \ref{lemma: mazurkiewicz extendability} we requested that the elements of $r$ are not in $\cup\tup{p_\beta}$. 
Notice that cardinality is not enough to argue this, since $p_\beta\supseteq p$ and $p$ can be of cardinality $\mathfrak{c}$.
This happens, for example, in the case that $p$ is a condition in $\mathbf{P_M}$.

\begin{Lema} \label{lemma: mazurkiewicz amalgamation}
	Let $\mathbf{Q}$ be the finite support product of $\omega_1$-many copies of Cohen forcing, let $g$ be $\mathbf{Q}$-generic over $V$.	
	Then $\mathbf{P_M}$ satisfies amalgamation in $V[g]$.
\end{Lema}

\begin{Dem}
	We need to prove that for densely many $p\in \mathbf{P_M}$, for any $g_1, g_2$ mutually $\mathbf{Q}$-generic over $V[p]$ and for all $p\in \mathbf{P_M}\cap V[p,g_1]$, $p_2\in \mathbf{P_M}\cap V[p,g_2]$ extending $p$, $p_1$ and $p_2$ are compatible.
	
	First, notice that $$D=\{ p \in \mathbf{P_M} \mid \exists \alpha<\omega_1 \, V[g\restriction\alpha] \models p \text{ is a Mazurkiewicz set} \}$$ is dense.
	For any condition $p\in \mathbf{P_M}$, there is a real $x$ such that $p\in V[x]$. 
	By Lemma \ref{lemma: prereq any real is in an inital segment of C(omega1)}, there is $\gamma<\omega_1$ such that $x\in V[g\restriction \gamma]$.
	Take $\bar{x}= \oplus_{\beta \leq \gamma} \cup g\restriction\{\beta\}$, and repeat the proof of Lemma \ref{lemma: mazurkiewicz extendability} for this $\bar{x}$. 
	Then there is $\bar{p}$ Mazurkiewicz set in $V[\bar{x}]=V[g\restriction \alpha]$ where $\alpha=\gamma + 1$. 
	So $\bar{p} \leq_{\mathbf{P_M}} p$ and $\bar{p}\in D$. 
	
	Now, fix $p\in D$ and let $x$ be such that $V[x]\models p $ is a Mazurkiewicz set.
	Let $g_1, g_2$ be mutually $\mathbf{Q}$-generic filters over $V[p]=V[x]$, and fix $p_1 \in \mathbf{P_M}\cap V[x,g_1]$ and  $p_2 \in \mathbf{P_M}\cap V[x,g_2]$ such that $p\subseteq p_1$, $p\subseteq p_2$.
	Let $y\in \RR\cap V[x, g_1]$ and $z\in \RR \cap V[x, g_2]$ be such that
	\begin{align*}
		V[x,y]\models& p_1 \text{ is a Mazurkiewicz set, and} \\
		V[x,z]\models& p_2 \text{ is a Mazurkiewicz set}.
	\end{align*}
	By Theorem \ref{th: prereq solovay}, we can choose $y$ and $z$ such that they are Cohen generic over $V[x]$. 
	Since $g_1$ and $g_2$ are mutually $\mathbf{Q}$-generic, $y$ and $z$ are mutually Cohen generic over $V[x]$. 
	
	We will show that $p_1\cup p_2$ is a partial condition in $V[x,y,z]$. 
	This is enough since, by Lemma \ref{lemma: mazurkiewicz extendability}, we can find a condition $\bar{p}$ that extends $p_1\cup p_2$, and therefore witnesses the compatibility between $p_1$ and $p_2$.
	
	Work in $V[x,y,z]$. 
	Suppose $p_1\cup p_2$ contains three different points on a line $l$. 
	Since $p_1$ and $p_2$ are partial conditions, each of these sets does not contain three collinear points. 
	Without loss of generality, we can assume $|l\cap p_2|=2$ and $|l\cap p_1|\geq 1$. 
	Notice that $|l\cap p_2|= 2$ implies that $l\in V[x,z]$. 
	Since $p_1\in V[x,y]$, we know that $|l\cap V[x,y]|\geq 1$.
	We divide in two cases, depending on whether $|l\cap V[x,y]|\geq 2$ or $|l\cap V[x,y]|=1$.

	\textbf{Case 1.}
	If $|l\cap V[x,y]|\geq 2$, then $l\in V[x,y]$, therefore $l\in V[x,y]\cap V[x,z]=V[x]$. 
	Since $p$ is a Mazurkiewicz set in $V[x]$, $|p\cap l|=2$. 
	Since $p_2\supseteq p$ is a Mazurkiewicz set in $V[x,y]$, $p_2\cap l = p\cap l$ and analogously for $p_1$.
	Therefore $(p_1\cup p_2) \cap l = p\cap l$, contradicting the choice of $l$. \\
	
	\textbf{Case 2.}	
	If $|l\cap V[x,y]|=1$, let $r$ be the only element in $l\cap V[x,y]$.
	Let $s_1, s_2\in l\cap p_2$.
	Then, 
	\[V[x,y,z] \models r \text{ is the only element of } V[x,y] \cap l(s_1,s_2). \]
	Recall that $y$ is generic over $V[x, z]$. There is a condition $t\in y \subseteq \mathbf{C}$ such that 
	\begin{equation} \label{eq: mazurkiewicz amalgamation}
		t \forc{\mathbf{C}}{V[x,z]} \dot{r}  \text{ is the only element of } V[x,\dot{g}] \cap l(\check{s}_1,\check{s}_2).
	\end{equation}
	
	Split $y$ in two mutually generic Cohen reals $y_1, y_2$ according to $t$ as in Definition \ref{def: prereq split of a cohen real}.
	From Equation \ref{eq: mazurkiewicz amalgamation}, we get that 
	\begin{align*}
		V[x,z,y_1] & \models r_1 \text{ is the only element of } V[x,y_1] \cap l(s_1,s_2), \text{ and} \\
		V[x,z,y_2] & \models r_2 \text{ is the only element of } V[x,y_2] \cap l(s_1,s_2);
	\end{align*}
	where $r_1= \dot{r}_{y_1}$ and $r_2= \dot{r}_{y_2}$.
	
	Since $V[x,y_1], V[x,y_2]\subseteq V[x,y]$ and $r, r_1, r_2 \in l$, we obtain that $r=r_1=r_2$. 
	Then, $r \in V[x,y_1]\cap V[x,y_2]$, so $r\in V[x]$.
	Thus, $(p_1\cup p_2 )\cap l = p_2\cap l$, contradicting the choice of $l$. \\
	
	Finally, there is no such $l$, and therefore $p_1\cup p_2$ is a partial condition. 
	By extendability (Lemma \ref{lemma: mazurkiewicz extendability}), there is $\bar{p}\in \mathbf{P_M}$ such that $\bar{p}\supseteq p_1\cup p_2 $. 
	Since the order in $\mathbf{P_M}$ is the reverse inclusion, we get $\bar{p}\leq_{\mathbf{P_M}} p_1$ and $\bar{p} \leq_{\mathbf{P_M}}p_2$. 
	Thus, $p_1$ and $p_2$ are compatible.
	
\end{Dem}

\begin{Coro}[of Theorem \ref{th: partition no wo}] \label{coro: mazurkiewicz - no wo}
	Let $\mathbf{Q}$ be the finite support product of $\omega_1$-many copies of  Cohen forcing, let $g$ be a $\mathbf{Q}$-generic filter over $V$. 
	Let $\mathbf{P}$ be the forcing poset in $V[g]$ given by
	\[p \in \mathbf{P} \iff \exists x \in \RR\; V[x]\models p \text{ is a Mazurkiewicz set}, \]
	ordered by reverse inclusion. 
	Let $h$ be a $\mathbf{P}$-generic filter over $V[g]$, and let $\mathcal{P}=\cup h$. 
	Then 
	\[L(\RR, \mathcal{P})^{V[g,h]}  \models \mathsf{ZF} + \mathsf{DC} +  \neg \mathsf{WO}(\RR) + \mathcal{P} \text{ is a Mazurkiewicz set}. \]
\end{Coro} 

\begin{Dem}
	We will use Theorem \ref{th: partition no wo}. 
	Notice that $\mathbf{P}=\mathbf{P_M}$, and we have shown that this forcing adds a real partition.
	Since the order in $\mathbf{P}$ is the reverse inclusion, we know that $\mathbf{P}$ is $\sigma$-closed applying Lemma \ref{lemma: sigma close for partitions}. 
	Finally, $\mathbf{P}$ satisfies extendability by Lemma \ref{lemma: mazurkiewicz extendability}, and amalgamation by Lemma \ref{lemma: mazurkiewicz amalgamation}.
	We can then apply Theorem \ref{th: partition no wo}, and obtain the desired conclusion. 
\end{Dem}

	\bibliographystyle{siam}
	\bibliography{../../bibliography.bib}
	
\end{document}